%% file: coordinates-3.tex
\documentclass[11pt]{article}
\usepackage{amsmath,amsfonts,epsfig}
\usepackage[all]{xy}
\usepackage{labelfig}
\usepackage{epstopdf}
\input{macros}

\widemargins

\newcommand{{\wln}}{wedgelation}
\newcommand{\x}{{\bf x}}
\newcommand{\y}{{\bf y}}

\begin{document}
\title{Shearing coordinates and convexity of length functions on Teichm\"uller space}

\author{M. Bestvina, K. Bromberg, K. Fujiwara, J. Souto\footnote{M.B. has been partially supported by NSF grant DMS-0502441, K.B. has been partially supported by NSF grants DMS-0554569 and DMS-0504877, K. F. has been partially supported by Grant-in-Aid for Scientific Research (No. 19340013) and J. S. has been partially supported by NSF grant DMS-0706878 and the Alfred P. Sloan Foundation.}}

\date{\today}

\maketitle

\begin{abstract}
\noindent We prove that there are Fenchel-Nielsen coordinates for the Teichm\"uller space of a finite area hyperbolic surface with respect to which the length functions are convex.
\end{abstract}

\section{Introduction}
Let $X$ be a complete hyperbolic surface with finite area and $\cT(X)$ its Teichm\"uller space. Recall that points in $\cT$ are equivalence classes of marked hyperbolic surfaces, i.e. equivalence classes of pairs $(Y,f)$ where $Y$ is a hyperbolic surface and $f:X\to Y$ is a quasi-conformal homeomorphism.

The goal of this paper is to describe certain coordinates of $\cT(X)$ with respect to which the length functions are convex. Given a homotopically essential, non-peripheral curve $\gamma\subset X$ and a point $(Y,f)\in\cT(S)$, then the curve $f(\gamma)$ is freely homotopic to a unique geodesic in the hyperbolic surface $Y$. Denoting by $\ell_\gamma(Y,f)$ the length of this geodesic we obtain a well-defined function
$$\ell_\gamma:\cT(X)\to\reals_+$$
This is the {\em length function} associated to the curve $\gamma$.

Of the many well-known ways to parametrize Teichm\"uller space perhaps the most classical is attributed to Fenchel and Nielsen. Fix a pants decomposition $\cP$ of $X$, i.e. a multicurve such that $X\setminus\cP$ is homeomorphic to the disjoint union of thrice punctured spheres. A simple computation shows that $\cP$ has $\vert\cP\vert=3g+n-3$ components where $g$ is the genus and $n$ the number of cusps of $X$. The Fenchel-Nielsen coordinates 
\begin{equation}\label{FN-coor}
\Phi:\cT(X)\xrightarrow{\sim}\reals_+^{3g+n-3}\times\reals^{3g+n-3}
\end{equation}
associate to each point in $\cT(X)$ the length and the twist for each component of the pants decomposition $\cP$. 

In the above discussion it the meaning of ''lenght'' is clear. To define the ''twist parameter''
it has to be decided what a $0$-twist represents. There is no canonical way to determine this so  
we refer to any of the possible length-twist coordinates as {\em Fenchel-Nielsen coordinates} associated to the pants decomposition $\cP$. We prove:

\begin{theorem}\label{thm:FN-convex}
Let $X$ be a complete, finite area, hyperbolic surface of genus $g$ and with $n$ cusps, and fix a pants decomposition $\cP$ of $X$. There are Fenchel-Nielsen coordinates $\Phi:\cT(X)\xrightarrow{\sim}\reals_+^{3g+n-3}\times\reals^{3g+n-3}$ associated to $\cP$ such that for any essential curve $\gamma$ in $X$ the function 
$$l_\gamma\circ\Phi^{-1}:\reals_+^{3g+n-3}\times\reals^{3g+n-3}\to\reals_+$$
is convex. If moreover the curve $\gamma$ intersects all the components of $\cP$ then $l_\gamma\circ\Phi$ is strictly convex.
\end{theorem}

Other convexity results of the length functions $l_\gamma$ are due to Kerckhoff \cite{Kerckhoff} and Wolpert \cite{Wolpert}. They proved respectively that the length functions are convex along earthquake paths and Weil-Petersson geodesics. Both authors derived from their results proofs of the so-called Nielsen realization problem; so do we.

\begin{theorem}[Kerckhoff]\label{realization}
The action of every finite subgroup of the mapping class group on $\cT(X)$ has a fixed point.
\end{theorem}

Tromba \cite{Tromba} gave a different proof of Theorem \ref{realization} using the convexity of the energy functional along Weil-Petersson geodesics. Proofs of this theorem in a completely different spirit are due to Gabai \cite{Gabai} and Casson-Jungreis \cite{Casson-Jungreis}.

In order to prove Theorem \ref{thm:FN-convex} we follow a slightly indirect path. We will associate a continuous map $s_\lambda:\cT(X)\to\reals^{\vert\lambda\vert}$ to every maximal, finite leaved lamination $\lambda$ of $X$. The image $T_\lambda=s_\lambda(\cT(X))$ of $s_\lambda$ is an open convex subset of a linear subspace of the correct dimension $6g+2n-6$. We refer to 
$$s_\lambda:\cT(X)\xrightarrow{\sim}T_\lambda$$ 
as {\em shearing coordinates} associated to the lamination $\lambda$. The shearing coordinates are closely related to those studied by Bonahon \cite{Bonahon}.

Given a pants decomposition $\cP$ of $X$, we choose a maximal lamination $\lambda$ containing $\cP$ and describe Fenchel-Nielsen coordinates associated to $\cP$ in such a way that the map 
$$\Phi\circ s_\lambda^{-1}:T_\lambda\xrightarrow{\sim}\reals_+^{3g+n-3}\times\reals^{3g+n-3}$$
is linear. In particular, Theorem \ref{thm:FN-convex} follows immediately from the following more general result:

\begin{theorem}\label{thm:shear-convex}
Let $X$ be a complete, finite area, hyperbolic surface with genus $g$ and with $n$ cusps. Let $\lambda$ be a maximal lamination in $X$ with finitely many leaves and let $s_\lambda:\cT(X)\xrightarrow{\sim} T_\lambda$ be the shearing coordinates associates to $\lambda$. For any essential curve $\gamma$ in $X$ the function 
$$l_\gamma\circ s_\lambda^{-1}:T_\lambda\to\reals_+$$
is convex. If moreover the curve $\gamma$ intersects all the leaves of $\lambda$ then $l_\gamma\circ s_\lambda^{-1}$ is strictly convex.
\end{theorem}

The paper is organized as follows. After a few preliminaries in section \ref{sec:preli}, we introduce in section \ref{sec:coordinates} the shearing coordinates and reduce the proof of Theorem \ref{thm:shear-convex} to Proposition \ref{lengthfactors}, our main technical result. In section \ref{sec:FN} we relate Fenchel-Nielsen coordinates to shearing coordinates and prove Theorem \ref{thm:FN-convex}. In section \ref{sec:annulus} we study the length function on the Teichm\"uller space of the annulus and finally in section \ref{sec:final} we prove Proposition \ref{lengthfactors}. The proofs are, once one is used to the notation, elementary.

%
%
%
%

\section{Teichm\"uller space and other important objects}\label{sec:preli}
Let $X$ be a complete, orientable hyperbolic surface with possibly infinite area. A {\em marked hyperbolic surface} is a pair $(Y,f)$ where $Y$ is a hyperbolic surface and $f: X \longrightarrow Y$ is a quasi-conformal homeomorphism. Two marked surfaces $(Y_0, f_0)$ and $(Y_1, f_1)$ are equivalent if there is an isometry $\phi :Y_0 \longrightarrow Y_1$ such that $\phi \circ f_0$ and $f_1$ are quasi-conformally isotopic. The Teichm\"uller space $\cT(X)$ is the set of equivalence classes of marked hyperbolic surfaces. We give $\cT(X)$ a metric (and topology) as follows. The distance, $d_T((Y_0,f_0), (Y_1, f_1))$, between two pairs is the infimum of the logarithm of the quasi-conformal constant of all maps $\phi: Y_0 \longrightarrow Y_1$ with $\phi \circ f_0$ quasi-conformally isotopic to $f_1$.
\medskip

\noindent{\bf Remark.} In this paper we will be mostly interested in Teichm\"uller spaces of finite area surfaces although the ambitious reader could easily see that our results can be extended to a more general setting. In particular, as we will see in section \ref{sec:annulus} the Teichm\"uller space of a hyperbolic annulus plays a key role in our work.
\medskip

Continuing with the same notation as above, we can identify the universal covers of $X$ and $Y$ with the hyperbolic plane $\htwo$; this identification is unique up to composition with an isometry of $\htwo$. Let $\tilde f:\htwo\to\htwo$ be the lift of the quasi-conformal homeomorphism $f:Y\to X$. It is well-known that $\tilde f$ extends continuously to a homeomorphism of $\del\tilde f:\del \htwo\to\del\htwo$; here $\del\htwo$ is the boundary at infinity of $\htwo$. Moreover, lifts of quasi-conformally isotopic maps have extensions which differ by composition with (the boundary extensions of) isometries of $\htwo$.


A {\em lamination} on $X$ is a closed, but perhaps not compact, subset of $X$ which is foliated by geodesics. If the surface $X$ has finite area, we will be only be interested in laminations $\lambda$ with only finitely many leaves. Recall that the only non-isolated leaves of such a lamination are simple closed geodesics in $X$.

Let $f: X \longrightarrow Y$ be a quasi-conformal homeomorphism between hyperbolic surfaces and let $\tilde{f}: \htwo \longrightarrow \htwo$ and $\del\tilde f:\del\htwo\to\del\htwo$ be as above. The maps $f$ and $\tilde{f}$ will not take geodesics to geodesics. In order to by-pass this problem, we associate to any geodesic $\gamma\subset\htwo$ the unique geodesic $\bar{f}(\gamma)\subset\htwo$ which has the same endpoints on $\del \htwo$ as the arc $\tilde{f}(\gamma)$. If $\gamma$ is a closed geodesic on $X$ then its pre-image $\tilde{\gamma}$ in $\htwo$ will be equivariant and therefore $\bar{f}(\tilde{\gamma})$ will also be equivariant and descend to a geodesic $\bar{f}(\gamma)$ on $Y$. If $\gamma$ is also simple then $\bar{f}(\gamma)$ is simple as well. For a lamination we apply $\bar{f}$ to each geodesic in the lamination. Similarly, for an ideal triangle $\Delta$ we let $\bar{f}(\Delta)$ be the ideal triangle whose boundary is the the $\bar{f}$-image of the geodesics bounding $\Delta$.

We refer to \cite{Imayosi-Taniguchi} and \cite{Casson-Bleiler} for facts and definitions in this section.

%
%
%
%

\section{Shearing coordinates}\label{sec:coordinates}
The goal of this section is to define {\em shearing coordinates} for the Teichm\"uller space $\cT(X)$ of a finite area hyperbolic surface.

\subsection{Ideal triangulations}
\begin{figure}[h]
        \centering
         \includegraphics[width=10cm]{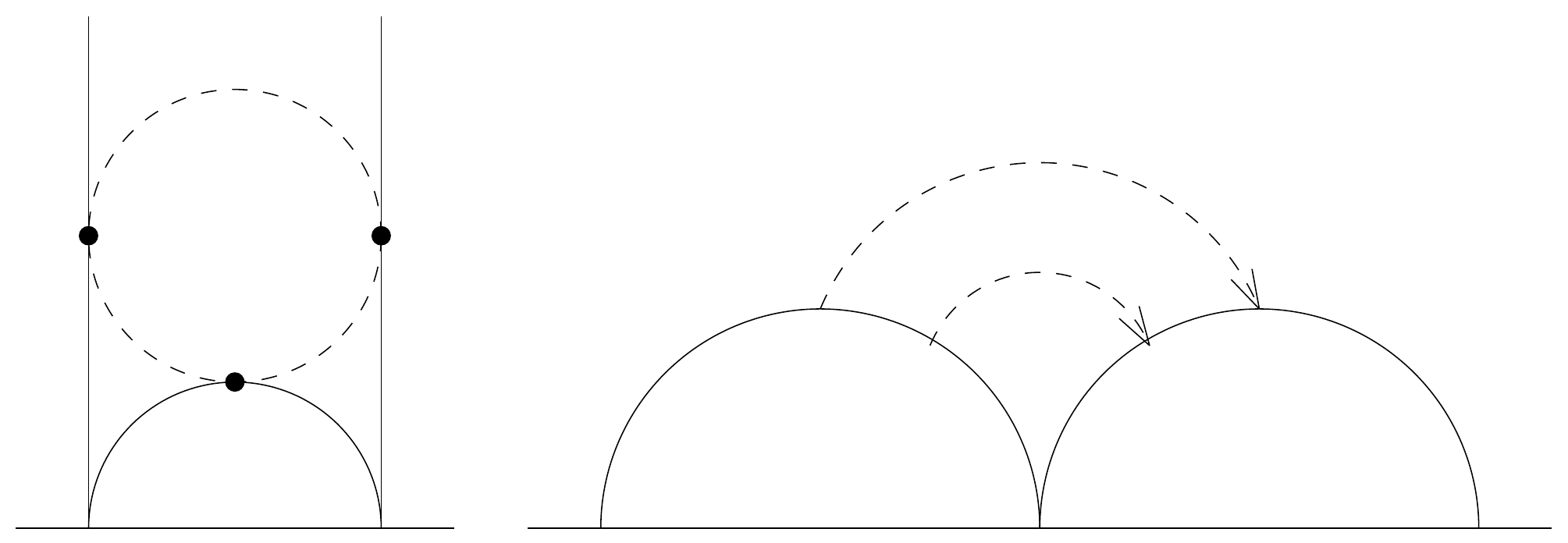}
         \caption{(a) and (b)}\label{fig1}
\end{figure}

Before setting up our coordinates we need some definitions and notation. An {\em ideal triangle} on a hyperbolic surface $X$ is the image of an injective, local isometry from an ideal triangle in $\htwo$ to the surface. Recall that any two ideal triangles in $\htwo$ are isometric and that ideal triangles have the same isometry group as euclidean equilateral triangles. An {\em ideal triangulation} of $X$ is a lamination with finitely many leaves whose complementary components are ideal triangles. An ideal triangle has a unique inscribed disk that is tangent to all three sides of the triangle. The {\em midpoints} of the sides are three tangency points; compare with figure \ref{fig1} (a). 

Let $\gamma_0$ and $\gamma_1$ be geodesics in $\htwo$ that are asymptotic to a point $p_\infty\in\del\htwo$. For each $p_0 \in \gamma_0$ there is a unique $p_1 \in \gamma_1$ such that the horocycle based at $p_\infty$ through $p_0$ intersects $\gamma_1$ at $p_1$. We define a map $h_{\gamma_0, \gamma_1}: \gamma_0 \longrightarrow \gamma_1$ by $h_{\gamma_0, \gamma_1}(p_0) = p_1$. We allow the possibility that $\gamma_0 = \gamma_1$ in which case $h_{\gamma_0, \gamma_1}$ is the identity map; compare with figure \ref{fig1} (b).

Let $\Delta^a$ and $\Delta^b$ be ideal triangles in $\htwo$ with disjoint interiors and let $\gamma$ be a geodesic separating the two triangles. We also assume that both $\Delta^a$ and $\Delta^b$ are asymptotic to $\gamma$; that is, there are sides $\gamma^a$ and $\gamma^b$ of $\Delta^a$ and $\Delta^b$ that are asymptotic to $\gamma$. Let $m^a$ and $m^b$ be the midpoints of $\gamma^a$ and $\gamma^b$.

\begin{figure}[h]
\leavevmode \SetLabels
\L (.51*.975) $\gamma$\\
\L (.35*.70) $\Delta^a$\\
\L (.61*.25) $\Delta^b$\\
\endSetLabels
\par
\begin{center}
\AffixLabels{\centerline{\epsfig{file =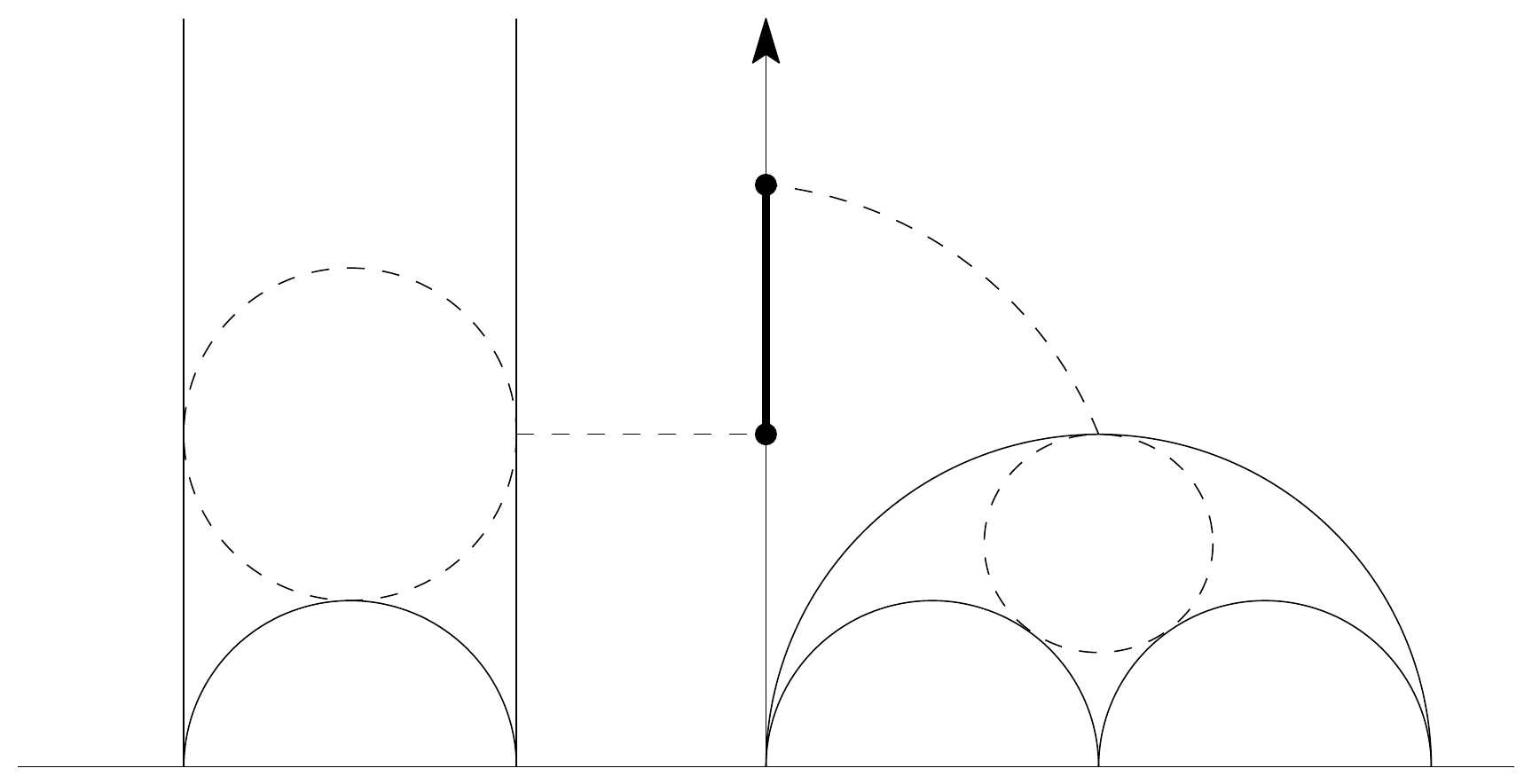, width=8cm, angle= 0}}}
\end{center}
\caption{The two dots are $h_{\gamma^a, \gamma}(m^a)$ and $h_{\gamma^b, \gamma}(m^b)$. In this picture $s(\Delta^a, \Delta^b, \gamma)$ is equal to the length of the bold printed segment.}\label{fig2}
\end{figure}

We define $s(\Delta^a, \Delta^b, \gamma)$ to be the signed distance between $h_{\gamma^a, \gamma}(m^a)$ and $h_{\gamma^b, \gamma}(m^b)$ where the sign is determined by orienting $\gamma$ such that $\Delta^a$ is on the left of $\gamma$; compare with figure \ref{fig2}. If $\Delta^a$ and $\Delta^b$ have a common boundary edge there is only one choice for $\gamma$ so we will sometimes write $s(\Delta^a, \Delta^b) = s(\Delta^a, \Delta^b, \gamma)$.


The following lemma is a collection of simple facts on $s(\cdot,\cdot,\cdot)$ whose proof we leave to the interested reader.

\begin{lemma}\label{lem:facts}
Let $\Delta^a$, $\Delta^b$ and $\gamma$ be as above.
\begin{itemize}
\item $s(\Delta^b, \Delta^a, \gamma) = s(\Delta^a, \Delta^b, \gamma)$.
\item If $\phi$ is an isometry of $\htwo$ then $s(\Delta^a, \Delta^b, \gamma) = s(\phi(\Delta^a), \phi(\Delta^b), \phi(\gamma)).$
\item $s(\Delta^a, \Delta^b, \gamma') = s(\Delta^a, \Delta^b, \gamma)$ for any another geodesic $\gamma'$ separating $\Delta^a$ and $\Delta^b$ and asymptotic to $\gamma$.
\end{itemize}
If moreover $\Delta^0,\Delta^1,\dots,\Delta^k$ is a chain of pairwise asymptotic ideal triangles with disjoint interior and $\gamma_i$ is for $i=1,\dots,k$ a geodesic separating $\Delta^0,\dots,\Delta^{i-1}$ from $\Delta^i,\dots,\Delta^k$ then we have
$$s(\Delta^0,\Delta^k,\gamma)=s(\Delta^0,\Delta^1,\gamma_1)+\dots+s(\Delta^{k-1},\Delta^k,\gamma_k)$$
for any geodesic $\gamma$ separating $\Delta^0$ and $\Delta^k$.\qed{lem:facts}
\end{lemma}

\begin{figure}[h]
\leavevmode \SetLabels
\L (.51*.655) $\gamma$\\
\L (.33*.70) $\Delta^a$\\
\L (.61*.25) $\Delta^b$\\
\endSetLabels
\par
\begin{center}
\AffixLabels{\centerline{\epsfig{file =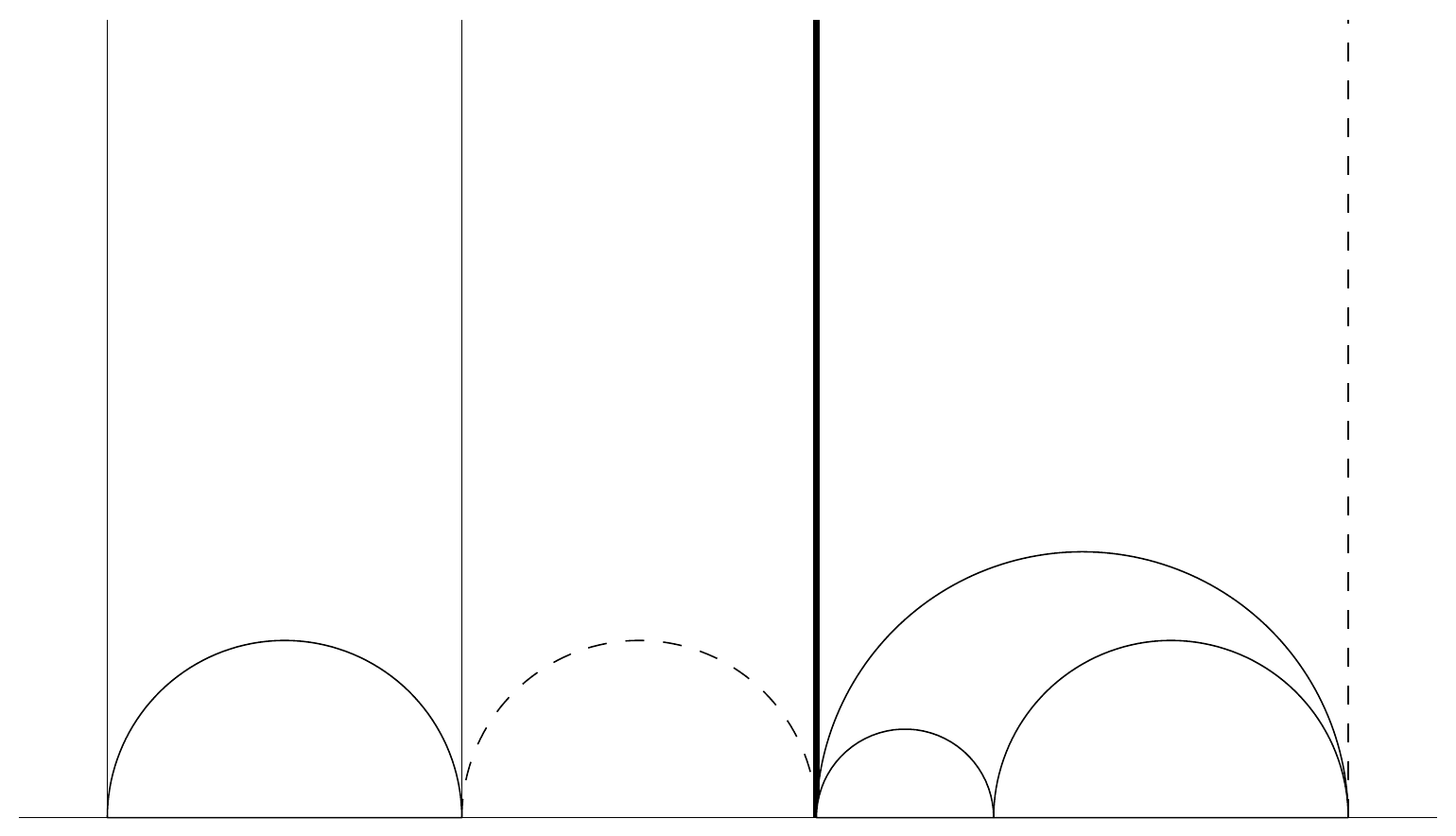, width=8cm, angle= 0}}}
\end{center}
\caption{The dotted lines determine the auxiliary ideal triangles needed to express $s(\Delta^a,\Delta^b,\gamma)$ in terms of cross ratios.}\label{fig7}
\end{figure}

\noindent{\bf Remark.} Working in say the upper half-plane model, assume that two ideal triangles $\Delta_a=[\theta_1,\theta_2,\theta_3]$ and $\Delta_b=[\theta_1,\theta_3,\theta_4]$ have a common boundary edge $[\theta_1,\theta_3]$. Let
$$\kappa=[\theta_1,\theta_2,\theta_3,\theta_4]=\frac{(\theta_1-\theta_3)(\theta_2-\theta_4)}{(\theta_1-\theta_4)(\theta_2-\theta_3)}$$
be the cross-ratio of the four vertices $\theta_1,\theta_2,\theta_3,\theta_4$. We have then the formula
$$s(\Delta^a,\Delta^b)=\log(\kappa-1)$$
Decomposing the general picture into adjacent ideal triangles and using the last claim of Lemma \ref{lem:facts}, it is easy to also express $s(\Delta^a,\Delta^b,\gamma)$ as sums of logarithms of algebraic expressions in cross ratios when $\Delta^a$ and $\Delta^b$ don't have a common boundary edge; compare with figure \ref{fig7}.


\subsection{The coordinates}
Let $\lambda$ be from now on an ideal triangulation of the complete finite area hyperbolic surface $X$. 
We now define a coordinate map $s_\lambda: \cT(X) \longrightarrow \reals^{|\lambda|}$ by defining a coordinate function $s_\gamma$ for each leaf $\gamma$ of $\lambda$. We need to make some choices in the definition of $s_\lambda$ but the map will be unique up to post-composition with a linear map.

Assume $\gamma$ is an isolated leaf of $\lambda$ and let $\tilde{\gamma}$ be a component of the pre-image of $\gamma$ in the universal cover. Let $\tilde{\lambda}$ be the pre-image of $\lambda$ in the universal cover. Let $\Delta^a_\gamma$ and $\Delta^b_\gamma$ be the two ideal triangle in the complement of $\tilde\lambda$ whose boundary contains $\tilde{\gamma}$. Then we define $s_\gamma(Y,f) = s(\bar{f}(\Delta^a_\gamma), \bar{f}(\Delta^b_\gamma), \bar{f}(\tilde{\gamma}))$. Here $\bar f$ is as in the end of section \ref{sec:preli}. It is clear that $s_\gamma(Y,f)$ is independent of the choice of the lift $\tilde{\gamma}$. 

Now assume that $\gamma$ is a closed curve in $\lambda$. To define $s_\gamma$ we need to make some arbitrary choices. Again let $\tilde{\gamma}$ be a component of the pre-image of $\gamma$ in the universal cover. In this case there will not be ideal triangles whose boundary contains $\tilde{\gamma}$. Instead we choose ideal triangles $\Delta^a_\gamma$ and $\Delta^b_\gamma$ such that $\tilde{\gamma}$ separates the triangles and they are both asymptotic to $\tilde{\gamma}$. We then define $s_\gamma(Y,f) = s(\bar{f}(\Delta^a_\gamma), \bar{f}(\Delta^b_\gamma), \bar{f}(\tilde{\gamma}))$. 
\medskip

\noindent {\bf Remark.} Observe that by the last claim of Lemma \ref{lem:facts}, different choices of $\Delta^a$ and $\Delta^b$, say $\hat\Delta^a$ and $\hat\Delta^b$, yield functions $\hat s_\gamma(Y,f) = s(\bar{f}(\hat\Delta^a_\gamma), \bar{f}(\hat\Delta^b_\gamma), \bar{f}(\tilde{\gamma}))$ which differ from $s_\gamma(Y,f)$ by a linear combination of the functions $s_\eta(Y,f)$ corresponding to the isolated leaves $\eta$ separating $\Delta^a$ from $\hat\Delta^a$ and $\Delta^b$ from $\hat\Delta^b$.
\medskip

We now define our coordinate map by
\begin{equation}\label{eq:shear-c}
s_\lambda: \cT(X) \longrightarrow \reals^{|\lambda|}
\end{equation}
by $s_{\lambda}(Y,f) = (s_{\gamma_1}(Y,f), \dots, s_{\gamma_n}(Y,f))$ where $\gamma_1, \dots, \gamma_n$ are the components of $\lambda$. By the remark above, different choices amount to postcomposing $s_\lambda$ with an invertible linear transformation of $\reals^{|\lambda|}$. 

To see that $s_\lambda$ is continuous we fix an identification of $\tilde{X} = \htwo$ with the upper half-plane. Given a quadruple $(\theta_1,\theta_2,\theta_2,\theta_4)$ with $\theta_i\in\del\htwo=\reals\cup\{\infty\}$ we then define a function
$$\cT(X)\to\reals,\ \ (Y,f)\mapsto [\del\tilde f(\theta_1),\del\tilde f(\theta_2),\del\tilde f(\theta_3),\del\tilde f(\theta_4)].$$
Since the cross-ratio is invariant under isometries of $\htwo$ this function is well-defined. It directly follows from the geometric definition of a quasi-conformal map that this function is continuous. Since
$s_\lambda(Y,f)$ can be expressed in terms of cross-ratios it follows that $s_\lambda$ is continuous. 
%

\subsection{Image of $s_\lambda$}
Our next goal is to determine the image of $s_\lambda$. Before doing so we need still some more notation. 

Assume that $\gamma$ is a non-isolated, i.e. closed, leaf of $\lambda$. Our choice of $\Delta^a_\gamma$ and $\Delta^b_\gamma$ determines an $a$-side and a $b$-side of $\gamma$. In particular if $C_\gamma$ is a collar neighborhood of $\gamma$ then $C_\gamma \backslash \gamma$ has two components which we call the {\em sides} of $\gamma$. If $\tilde{C}_\gamma$ is the component of the pre-image of $C_\gamma$ that contains $\tilde{\gamma}$ then $\Delta^a_\gamma$ will intersect one of the components of $\tilde{C}_\gamma  \backslash \tilde\gamma$. The image of this component in $X$ will be one of the components of $C_\gamma \backslash \gamma$. This is the $a$-side. Then $\Delta^b_\gamma$ will intersect the other component of $\tilde{C}_\gamma \backslash \tilde{\gamma}$ and this component will map to the $b$-side of $\gamma$.

We need to assign a sign to each side of $\gamma$ that will be determined by the direction $\Delta^a_\gamma$ and $\Delta^b_\gamma$ spiral around $\gamma$. As above, orient $\tilde{\gamma}$ so that $\Delta^a_\gamma$ is on the left. Then $\sigma^a_\gamma = -1$ if $\Delta^a_\gamma$ is asymptotic to the forward end of $\tilde{\gamma}$ and $\sigma^a_\gamma = +1$ if $\Delta^a_\gamma$ is asymptotic to the negative end. We make a similar definition for $\sigma^b_\gamma$. Note that since $\tilde{\gamma}$ separates $\Delta^a_\gamma$ form $\Delta^b_\gamma$ we need to change the orientation of $\tilde\gamma$ when we define $\sigma^b_\gamma$.


\begin{figure}[h]
\leavevmode \SetLabels
\L (.50*.655) $\gamma$\\
\L (.35*.40) $\Delta^a$\\
\L (.63*.40) $\Delta^b$\\
\endSetLabels
\par
\begin{center}
\AffixLabels{\centerline{\epsfig{file =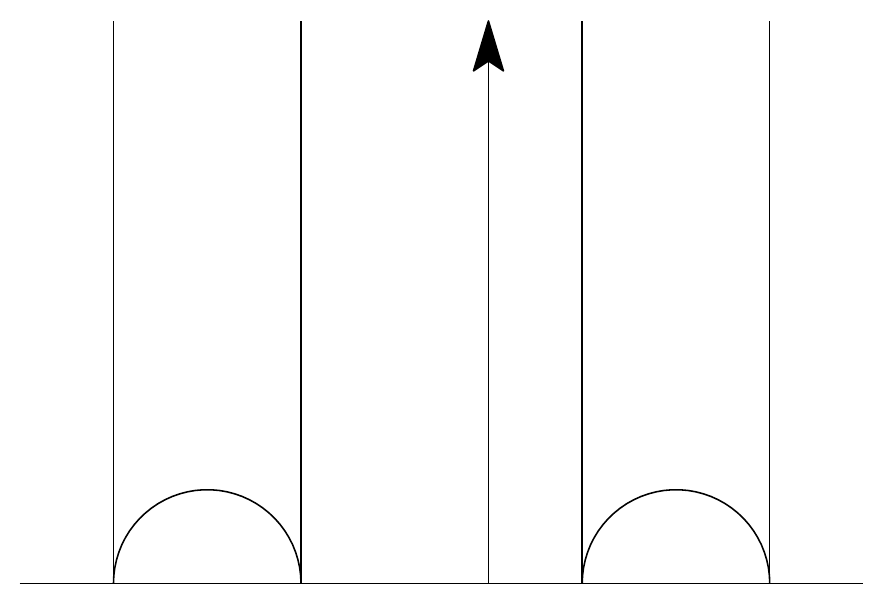, width=8cm, angle= 0}}}
\end{center}
\caption{In this picture $\sigma^a_\gamma=-1$ and $\sigma^b_\gamma=1$.}\label{fig5}
\end{figure}

Still assuming that $\gamma$ is a non-isolated leaf of $\lambda$, we will now express its length of in terms of our coordinates. The intersection of each leaf $\beta$ of $\lambda$ with the $a$-side of $\gamma$ will have $0$, $1$ or $2$ components. Let $n^a_\gamma(\beta)$ be this number and similarly define $n^b_\gamma(\beta)$. Note that if $\beta$ is a closed leaf then the number will always be zero. The content of the following lemma is that the length function $\ell_\gamma$ is given by any of the following two linear functions
$$\ell^a_\gamma:\reals^{|\lambda|}\to\reals,\ \ \ell^a_\gamma(\x) = \sigma^a_\gamma \sum_{i=1}^{|\lambda|} n^a_\gamma(\gamma_i) x_i$$
$$\ell^b_\gamma:\reals^{|\lambda|}\to\reals,\ \ \ell^b_\gamma(\x) = \sigma^b_\gamma \sum_{i=1}^{|\lambda|} n^b_\gamma(\gamma_i) x_i$$

\begin{lemma}\label{gammalength}
If $\gamma$ is a closed curve in $\lambda$ and $\ell_\gamma: \cT(X) \longrightarrow \reals$ is its length function then
$$\ell_\gamma = \ell^a_\gamma \circ s_\lambda = \ell^b_\gamma \circ s_\lambda.$$
\end{lemma}

\noindent {\bf Proof.} Let $(Y,f)$ be a marked hyperbolic structure in $\cT(X)$. We will calculate $\ell^a_\gamma \circ s_\lambda(Y,f)$. Working in the universal cover $\htwo$ using the upper half space model we can assume that $\bar{f}(\tilde{\gamma})$ is the vertical line at $x = 0$.

If we assume that $\sigma^a_\gamma = -1$ then we can also choose $\tilde{f}$ so that the geodesics in $\bar{f}(\tilde{\lambda})$ that intersect that $a$-side of $\bar{f}(\tilde{\gamma})$ are vertical geodesics with negative $x$-coordinates. Label these $x$-coordinates $x_i$ with $x_{i+1}<x_i$. By our normalization the subgroup of the deck group for $Y$ that fixes $\bar{f}(\tilde{\gamma})$ will be generated by the isometry $z \mapsto e^{\ell_\gamma(Y,f)} z$. The set of vertical geodesics will be invariant under this isometry and we will have $x_{i+k} = e^{\ell_\gamma(Y,f)} x_i$ where
$$k = \sum_{\beta \in \lambda} n^a_\gamma(\beta)$$
is the number of components of the intersection of $\lambda$ with the $a$-side of $\gamma$. Therefore $\ell_\gamma(Y,f) = \log x_{i+k}/x_i$.

Let $\Delta_i$ be the ideal triangle that has two vertical sides with $x$-coordinate $x_i$ and $x_{i+1}$. Note that the midpoints of the two vertical sides will have the same $y$-coordinate which we label $m_i$ and that $x_{i+1} = x_i - m_i$. We also observe that $s(\Delta_i, \Delta_{i+1}) = -\log m_{i+1}/m_i$ so
\begin{eqnarray*}
\sum_{i=0}^{k-1} s(\Delta_i, \Delta_{i+1}) & = & -\sum_{i=0}^{k-1} \log m_{i+1}/m_i \\
& = & -\log m_{k} + \log m_0 \\
& = &  -\log (x_k - x_{k+1}) + \log (x_0 - x_1)\\
& = & -\log x_k/x_0 \\
&=& -\ell_\gamma(Y,f).
\end{eqnarray*}

The second to last equality follows from the fact that $\ell_\gamma(Y,f) = \log x_k/x_0 = \log x_{k+1}/x_1$ and therefore $x_{k+1}/x_k = x_1/x_0$.

To finish the proof we need to write the sum on the left in terms of our coordinates. To do so we note that each vertical geodesic with $x$-coordinate $x_i$ maps to an isolated component $\beta$ of $\lambda$ and $s_\beta(Y,f) = s(\Delta_{i-1}, \Delta_i)$. Furthermore this geodesic will intersect the $a$-side of $\gamma$ so $n^a_\gamma(\beta)$ is positive. In fact each component $\beta$ of $\lambda$ that intersects the $a$-side of $\gamma$ will have exactly $n^a_\gamma(\beta)$ pre-images among the vertical geodesics with $x$-coordinates $x_1, \dots, x_k$. Therefore
\begin{eqnarray*}
\ell^a_\gamma \circ s_\lambda (Y,f) & = & \sigma^a_\gamma \sum_{i=1}^{|\lambda|} n^a_\gamma(\gamma_i) s_{\gamma_i}(Y,f) \\
& = & -\sum_{\beta \in \lambda} n^a_\gamma(\beta) s_\beta(Y,f) \\
& = & -\sum_{i=0}^{k-1} s(\Delta_i, \Delta_{i+1})
\end{eqnarray*}
which completes the proof when $\sigma^a_\gamma = -1$.

When $\sigma^a_\gamma = +1$ the proof is exactly the same except the vertical geodesics on the $a$-side of $\bar{f}(\tilde{\gamma})$ have positive $x$-coordinate. If we label the coordinates $x_i$ with $x_i <x_{i+1}$ then $s(\Delta_i, \Delta_{i+1}) = \log m_{i+1}/m_i$. The rest of the proof of the proof is exactly the same so this accounts for the $\sigma^a_\gamma= +1$ in the definition of $\ell^a_\gamma$.

The proof for $\ell^b_\gamma$ is obviously the same. \qed{gammalength}

Now let $c$ be a cusp of $X$. We let $n_c(\beta)$ be the number of components of the intersection of $\beta$ with a horospherical neighborhood of $c$ and define
$$\ell_c(\x) = \sum_{i=1}^{|\lambda|} n_c(\gamma_i) x_i.$$
We then have the following lemma.
\begin{lemma}\label{cusplength}
If $c$ is a cusp of $X$ then $\ell_c \circ s_\lambda(Y,f) = 0$ for all $(Y,f)$ in $\cT(X)$.
\end{lemma}

\noindent {\bf Proof.} The proof follows the same basic idea as the proof of Lemma \ref{gammalength}. We can assume that a component of the pre-image of the horosphere neighborhood of the cusp is a horosphere neighborhood of infinity in the upper half-space model of $\htwo$. Then the geodesics in $\bar{f}(\tilde{\lambda})$ that intersect this neighborhood will be vertical geodesics with $x$-coordinates $x_i$ labeled such that $x_i < x_{i+1}$. The neighborhood will be invariant under the deck transformation that fixes infinity and we can assume that it is of the form $z \mapsto z +1$. In particular $x_{i+k} = x_i + 1$ where
$$k = \sum_{\beta \in \lambda} n_c(\beta).$$
The set of midpoints $m_i$ will also be invariant under the action of $z \mapsto z+1$ so we also have $m_i = m_{i+k}$. Repeating the calculations from Lemma \ref{gammalength} we see that
$$\sum_{i=0}^{k-1} s(\Delta_i,\Delta_{i+1}) = \log x_k/x_0 = 0$$
and that
$$\ell_c \circ s_\lambda(Y,f) = \sum_{i=0}^{k-1} s(\Delta_i, \Delta_{i+1})$$
and the lemma is proven. \qed{cusplength}

The functions $\ell^a_\gamma, \ell^b_\gamma$ and $\ell_c$ are linear functions defined on all of $\reals^{|\lambda|}$, in particular we deduce from Lemma \ref{gammalength} and Lemma \ref{cusplength}:

\begin{cor}\label{cor:image}
The image of $s_\lambda$ is contained in the convex polytope $T_\lambda\subset\reals^{|\lambda|}$ consisting of all those $\x$ with $\ell^a_\gamma(\x) = \ell^b_\gamma(\x)>0$ for all simple closed curves $\gamma$ in $\lambda$ and $\ell_c(\x) = 0$ for all cusps of $X$.\qed{cor:image}
\end{cor}


\subsection{The coordinates are coordinates}
We prove now the central result of this section:

\begin{theorem}\label{homeomorphism}
The map $s_\lambda$ is a homeomorphism from $\cT(X)$ to $T_\lambda$.
\end{theorem}

We will derive Theorem \ref{homeomorphism} from some standard results about Teichm\"uller space and the following proposition.

\begin{prop}\label{lengthfactors}
Let $\alpha$ be a closed curve on $X$ and $\ell_\alpha : \cT(X) \longrightarrow \reals^+$ its length function. Then there is a convex function $\bar{\ell}_\alpha: \reals^{|\lambda|} \longrightarrow \bar{\reals}$ such that $\ell_\alpha = \bar{\ell}_{\alpha} \circ s_\lambda$. Furthermore the $\bar{\ell}_{\alpha}$-image of $T_\lambda$ is in $\reals^+$ and if $\alpha$ intersects every leaf of $\lambda$ then $\bar{\ell}_\alpha$ is strictly convex.
\end{prop}

Proposition \ref{lengthfactors} is the main work of this paper; we defer its proof to the final section.
\medskip

\noindent {\bf Proof of Theorem \ref{homeomorphism}.} As mentioned before, the fact that $s_\lambda$ can be expressed in terms of cross ratios shows that it is continuous. By Corollary \ref{cor:image} the image lies in $T_\lambda$ which is a convex subset of a linear subspace of $\reals^{|\lambda|}$. The number of ideal triangles in $X \backslash \lambda$ is $-2\chi(X) = \area(X)/\pi$ and therefore $|\lambda| = -3\chi(X) + n_c$ where $n_c$ is the number of closed curves in $\lambda$. In the definition of $T_\lambda$ there is one equation for each closed curve and one equation for each cusp. These equations are also linearly independent so $T_\lambda$ is an open subset of a linear subspace of dimension $|\lambda| - n_c - n_p$ where $n_p$ is the number of punctures (or cusps). Note that $\dim \cT(X) = -3\chi(X) - n_p = |\lambda| - n_c - n_p$.

If $s_\lambda(Y_0) = s_\lambda(Y_1)$ then by Proposition \ref{lengthfactors}, $\ell_\alpha(Y_0) = \ell_\alpha(Y_1)$ for all closed curves $\alpha$ and therefore $Y_0 = Y_1$. This implies that $s_\lambda$ is injective. By invariance of domain $s_\lambda$ is a local homeomorphism and its image an open subset of $T_\lambda$.

To see that the image is closed take a finite collection of closed curves $\Gamma = \{\alpha_1, \dots, \alpha_k\}$ such that the complementary pieces are disks or punctured disks. Such a collection {\em binds} the surface. Let $\ell_\Gamma: \cT(X) \longrightarrow \reals^+$ be the sum of the length functions $\ell_{\alpha_i}$ and similarly define $\bar{\ell}_\Gamma:\reals^{|\lambda|} \longrightarrow \reals$ and the sum of the $\bar{\ell}_{\alpha_i}$. 

Let $Y_i$ be a sequence in $\cT(X)$ such that $s_\lambda(Y_i)$ converges to a point $\x$ in $T_\lambda$. By Proposition \ref{lengthfactors} $\bar{\ell}_\Gamma(\x) < \infty$ and since $\ell_\Gamma(Y_i) = \bar{\ell}_\Gamma \circ s_\lambda(Y_i) \to \bar{\ell}_\Gamma(\x)$ we see that $\ell_\Gamma(Y_i)$ is uniformly bounded. Therefore the $Y_i$ lie in a compact subset of $\cT(X)$ by Lemma 3.1 of \cite{Kerckhoff} or Proposition 2.4 of \cite{Thurston}. In particular we can find a subsequence such that $Y_{i_k}$ converges to some $Y \in \cT(X)$ and therefore $\lim s_\lambda(Y_i) = \lim s_\lambda(Y_{i_k}) = s_\lambda(Y)$ is in image of $s_\lambda$. Therefore the image is a closed set.

Since the image is open and closed it must be all of $T_\lambda$. \qed{homeomorphism}

Before moving observe that Theorem \ref{thm:shear-convex} in an immediate consequence of Proposition \ref{lengthfactors} and Theorem \ref{homeomorphism}:
\medskip

\noindent {\bf Theorem \ref{thm:shear-convex}} {\em 
Let $X$ be a complete hyperbolic surface with finite area of genus $g$ and with $n$ cusps. Let $\lambda$ be a maximal lamination in $X$ with finitely many leaves and let $s_\lambda:\cT(X)\xrightarrow{\sim} T_\lambda$ be the shearing coordinates associates to $\lambda$. For any essential curve $\gamma$ in $X$ the function 
$$l_\gamma\circ s_\lambda^{-1}:T_\lambda\to\reals_+$$
is convex. If moreover the curve $\gamma$ intersects all the leaves of $\lambda$ then $l_\gamma\circ s_\lambda^{-1}$ is strictly convex.\qed{thm:shear-convex}}

\subsection{The Nielsen realization conjecture}
Kerckhoff proved the follow theorem using the convexity of length functions along earthquake paths. We give a similar proof using our convexity result.
\begin{theorem}[Kerckhoff]\label{uniquemin}
If $\Gamma$ binds then $\ell_\Gamma$ has a unique minimum on $\cT(X)$.
\end{theorem}

\noindent {\bf Proof.} Fix an ideal triangulation $\lambda$. Since $\Gamma$ binds there must be some curve in $\Gamma$ that is not in $\lambda$ so by Proposition \ref{lengthfactors} the length function $\bar{\ell}_\Gamma$ is strictly convex on $T_\lambda$. As we noted in the proof of Theorem \ref{homeomorphism} the function $\ell_\Gamma$ is proper and since $\ell_\Gamma = \bar{\ell}_\Gamma \circ s_\lambda$ the function $\bar{\ell}_\Gamma$ is also proper. A proper, strictly convex function that is bounded below has a unique minimum so $\bar{\ell}_\Gamma$, and therefore $\ell_\Gamma$, have unique minimums. \qed{uniquemin}

Kerckhoff used this result to prove the Nielsen realization conjecture. As the proof is short we include it here.
\medskip

\noindent {\bf Theorem \ref{realization} (Kerckhoff)}
{\em The action of every finite subgroup of the mapping class group on $\cT(S)$ has a fixed point.}
\medskip

\noindent {\bf Proof.} Let $G$ be the finite subgroup. The $G$ orbit of any finite set of binding curves will still bind and will be $G$-invariant. Let $\Gamma$ be such a $G$-invariant binding set and let $X$ be the unique minimum of $\ell_\Gamma$. Clearly $X$ is fixed by $G$. \qed{realization}

\subsection{Surfaces with geodesic boundary}
We conclude this section with a few remarks on surfaces with boundary. Let $X$ be a finite area hyperbolic surface with boundary components $\beta_1, \dots, \beta_k$. For each boundary component $\beta_i$ choose an interval $I_i$ of $(0,\infty)$ and let $\cT(X; I_1, \dots, I_k)$ be the Teichm\"uller space of marked hyperbolic surfaces where the length of $\beta_i$ is in the interval $I_i$. We allow the $I_i$ to be open, closed, half-open or a point.

An ideal triangulation of $X$ is still a finite leaved geodesic lamination whose complement is the union of an open ideal triangle; in particular, $\partial X\subset\lambda$. Let $\hat\lambda=\lambda\setminus\partial X$ be the union of the all interior leaves of $\lambda$. The definition of the coordinate map $s_\lambda: \cT(X; I_1, \dots, I_k) \longrightarrow \reals^{|\hat\lambda|}$ still makes sense. For each boundary component $\beta$ of $X$ we have a function $\ell_\beta: \reals^{|\hat\lambda|} \longrightarrow \reals$ where the composition $\ell_\beta \circ s_\lambda$ is the length function for the boundary component. The definition of $\ell_\beta$ is exactly the same as the definition of the functions $\ell^a_\gamma$ and $\ell^b_\gamma$. Let $T_{\lambda,I_1,\dots,I_k}$ be the subset of $T_\lambda\subset\reals^{|\hat\lambda|}$ satisfying $\ell_{\beta_i}(\x) \in I_i$ for $i=1,\dots,k$.
 
\begin{theorem}\label{boundary}
The map $s_\lambda$ is a homeomorphism from $\cT(X; I_1, \dots, I_k)$ to $T_{\lambda,I_1,\dots,I_k}$.
\end{theorem}

\noindent {\bf Proof.} Note that if each $I_i$ is a point then $T_\lambda$ is an open subset of a linear subspace of $\reals^{|\hat\lambda|}$ and the proof is the same as the proof Theorem \ref{homeomorphism}. The general case follows from this observation. \qed{boundary}

%
%
%
%

\section{Fenchel-Nielsen coordinates}\label{sec:FN}
In this section we prove Theorem \ref{thm:FN-convex}. We will assume knowledge of some basic form of the Fenchel-Nielsen coordinates; see for example \cite{Buser}.
\medskip

\begin{figure}[h]
\begin{center}
\AffixLabels{\centerline{\epsfig{file =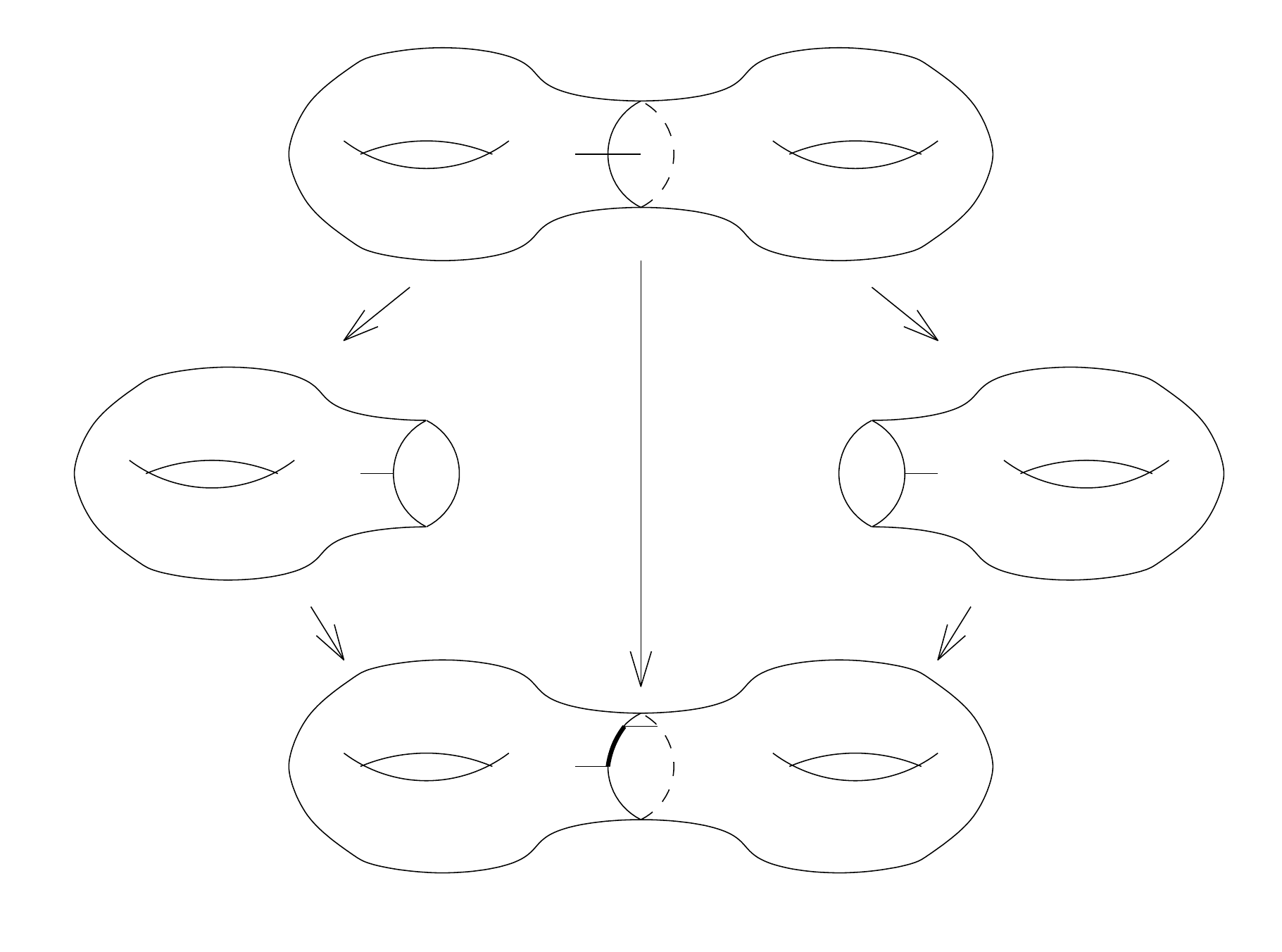, width=6cm, angle= 0}}}
\end{center}
\caption{Twist deformation for $t$ equal to the length of the bold printed arc}\label{fig8}
\end{figure}

Let $X$ be a hyperbolic surface and $\gamma$ a non-peripheral simple closed curve. Fenchel and Nielsen defined for $t\in\reals$ the twist deformation $T_\gamma^t(X)$ of $X$ along $\gamma$ with twist parameter $t$ as follows. First isotope $\gamma$ to the geodesic in its free homotopy class, also denoted by $\gamma$, and let $g:\reals\to\gamma$ be a parametrization by arc-length. Let $\gamma_1$ and $\gamma_2$ be the two boundary components obtained after cutting $X$ along $\gamma$ and for $i=1,2$ let $g_i:\reals\to\gamma_i$ be the induced parametrization. Up to relabeling, we may assume that $g_1$ is orientation preserving and $g_2$ is orientation reversing with respect to the induced orientation of the boundary curves $\gamma_1$ and $\gamma_2$. The hyperbolic surface $T_\gamma^t(X)$ is obtained from the cut open surface by identifying the points $g_1(s)$ with $g_2(s+t)$ fro all $s \in \reals$. Compare with figure \ref{fig8}.

The surface $T^t_\gamma(X)$ is just a hyperbolic surface; in particular it has so far no marking. However, it is well-known that there is flow, the Fenchel-Nielsen Dehn-twist flow
$$\tau_\gamma:\reals\times\cT(X)\to\cT(X),\ \ (t,(Y,f))\mapsto\tau_\gamma^t(Y,f)$$
such that for $(Y,f)\in\cT(X)$, the surface $T_\gamma^t(Y)$ is the hyperbolic surface associated to the point $\tau_\gamma^t(Y,f)$.

If $\gamma$ and $\gamma'$ are disjoint curves, then the flows $\tau_\gamma$ and $\tau_{\gamma'}$ commute. In particular, labeling by $\gamma_1,\dots,\gamma_{\vert\cP\vert}$ the components of some pants decomposition $\cP$ of $X$ we have an action
$$\tau_\cP:\reals^{|\cP|}\times\cT(X)\to\cT(X),\ \ ((t_1,\dots,t_{|\cP|}),(Y,f))\mapsto\tau_{\cP}^{(t_1,\dots,t_{|\cP|})}(Y,f)$$
where $\tau_{\cP}^{(t_1,\dots,t_{|\cP|})}=\tau_{\gamma_1}^{t_1}\circ\dots\circ\tau_{\gamma_{\vert\cP\vert}}^{t_{\vert\cP\vert}}$. It is well-known that the action given $\tau_{\cP}$ is free and proper.

To the pants decomposition $\cP$ of $X$ we can also associate the function
$$\ell_\cP:\cT(X) \longrightarrow (\reals^+)^{|\cP|}$$
which assigns to each point in $\cT(X)$ the $|\cP|$-tuple of lengths of curves in the pants decomposition.  Observe that the $\reals^{\vert\cP\vert}$-action $\tau_{\cP}$ preserves by definition the function $\ell_{\cP}$ and hence acts on the fibers. 

In fact, any existence results for Fenchel-Nielsen coordinates implies that the $\tau_{\cP}$-orbits actually coincide with the fibers of $\ell_{\cP}$. We can summarize this discussion as follows:

\begin{prop}\label{prop:bundle}
Let $\cP$ be a pants decomposition of a finite area surface $X$. Then
$$\ell_{\cP}:\cT(X)\to(\reals^+)^{\vert\cP\vert}$$
has a natural structure as a $\reals^{\vert\cP\vert}$-principal bundle.
\end{prop}

Every principal bundle over a contractible space is trivial but not canonically trivialized. In fact, any choice of a section yields a trivialization and vice-versa. From the point of view of the authors, all coordinates for Teichm\"uller space obtained by trivializing the principal bundle $\ell_{\cP}:\cT(X)\to(\reals^+)^{\vert\cP\vert}$ deserve to be referred to as {\em Fenchel-Nielsen} coordinates. We are now ready to prove Theorem \ref{thm:FN-convex}.
\medskip

\noindent{\bf Theorem \ref{thm:FN-convex}.}
{\em Let $X$ be a complete, finite area, hyperbolic surface of genus $g$ and with $n$ cusps, and fix a pants decomposition $\cP$ of $X$. There are Fenchel-Nielsen coordinates $\Phi:\cT(X)\xrightarrow{\sim}\reals_+^{3g+n-3}\times\reals^{3g+n-3}$ associated to $\cP$ such that for any essential curve $\gamma$ in $X$ the function 
$$l_\gamma\circ\Phi^{-1}:\reals_+^{3g+n-3}\times\reals^{3g+n-3}\to\reals_+$$
is convex. If moreover the curve $\gamma$ intersects all the components of $\cP$ then $l_\gamma\circ\Phi$ is strictly convex.}
\medskip

To begin with extend $\cP$ to an ideal triangulation $\lambda$ and let 
$$s_\lambda:\cT(X)\to T_\lambda$$ 
be the shearing coordinates associated to $\lambda$. Theorem \ref{thm:FN-convex} will follow immediately from Theorem \ref{thm:shear-convex} when we interpret the shearing coordinates $s_\lambda$ as Fenchel-Nielsen coordinates. 

Denote by $\lambda_0=\lambda\setminus\cP$ the set of isolated leaves in $\lambda$, recall the definition of the convex polytope $T_\lambda\subset\reals^{\vert\lambda\vert}=\reals^{\vert\cP\vert}\times\reals^{\vert\lambda_0\vert}$ and observe that the factor $\reals^{\vert\cP\vert}\times\{0\}$ in the above splitting of $\reals^{\lambda}$ is contained in $T_\lambda$. In particular, the canonical $\reals^{\vert\cP\vert}$-principal bundle structure on 
$$\pi:\reals^{\vert\lambda\vert}=\reals^{\vert\cP\vert}\times\reals^{\vert\lambda_0\vert}\to\reals^{\vert\lambda_0\vert}$$
induces a $\reals^{\vert\cP\vert}$-principal bundle structure on 
$$\pi:T_\lambda\to\pi(T_{\lambda})$$
Here $\pi$ is the projection to the second factor of the splitting of $\reals^{\vert\lambda\vert}$.

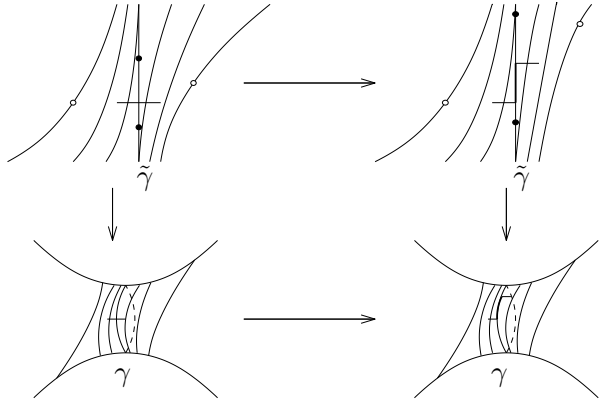
\begin{figure}[h]
\leavevmode \SetLabels
\L (.33*.055) $\gamma$\\
\L (.67*.055) $\gamma$\\
\L (.35*.54) $\tilde\gamma$\\
\L (.69*.54) $\tilde\gamma$\\
\endSetLabels
\par
\begin{center}
\AffixLabels{\centerline{\epsfig{file =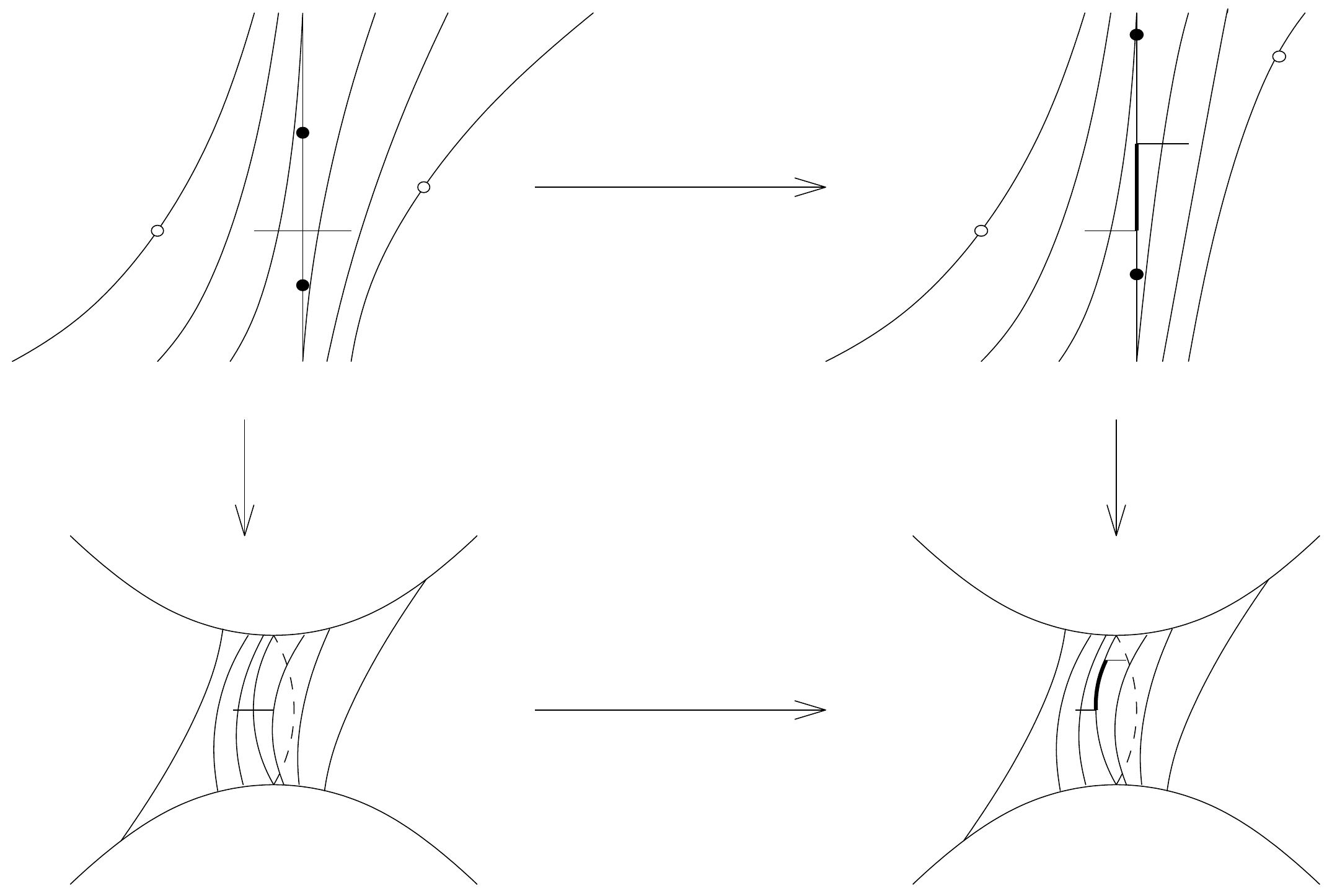, width=8cm, angle= 0}}}
\end{center}
\caption{The lower horizontal arrow is the twist deformation along the central curve $\gamma$ where the parameter is the length of the bold printed arc. The pictures in the upper row schematize the situation in the universal cover: the straight lines represent a lift $\tilde\gamma$ of $\gamma$, the white dots are the midpoints of the sides of the chosen triangles $\Delta^a$ and $\Delta^b$; the black dots are the projections of these midpoints to $\tilde\gamma$. The increase in the shearing coordinate associated to $\gamma$ coincides with the twist parameter.}\label{fig9}
\end{figure}

It follows directly from the definition just before Lemma \ref{gammalength} that for all $\gamma\in\cP$ the linear form $\ell_\gamma^a:T_\lambda\to\reals^+$ is independent of the $\reals^{\vert\cP\vert}$-factor and hence induces a well-defined linear function $\hat\ell_\gamma$ on $\pi(T_\lambda)$. Observe that this function is positive by Lemma \ref{gammalength}. Denoting by $\hat\ell_{\cP}:\pi(T_\lambda)\to(\reals^+)^{\vert\cP\vert}$ the function whose $\gamma$-coordinate is $\hat\ell_\gamma$ we have from Lemma \ref{gammalength} that the following diagram commutes
\begin{equation}\label{eq:bundles}
\xymatrix{
\cT(X)\ar[r]^{s_\lambda}\ar[d]_{\ell_{\cP}} & T_\lambda\ar[d]^\pi\\
(\reals^+)^{\vert\cP\vert} & \pi(T_\lambda)\ar[l]_{\hat\ell_{\cP}}
 }
\end{equation}
The map $\hat\ell_{\cP}$ is then linear and surjective. Since
$$\dim(\pi(T_\lambda))\le\dim(T_\lambda)-\dim(\reals^{\vert\cP\vert})\le\vert\cP\vert$$
we obtain that $\hat\ell_{\cP}$ is a linear isomorphism. In particular, the bundles $\cT(X)\to(\reals^+)^{\vert\cP\vert}$ and $T_\lambda\to\pi(T_{\lambda})$ are isomorphic as fiber bundles. Moreover, it follows directly from the definition of the shearing coordinates that they are also equivalent as $\reals^{\vert\cP\vert}$-principal bundles. Compare with figure \ref{fig9}.

The projection $T_\lambda\to(\reals^+)^{\vert\cP\vert}$ is linear and has kernel of dimension $\vert\cP\vert$. In particular, there is a linear map $L:T_\lambda\to\reals^{\vert\cP\vert}\times(\reals^+)^{\vert\cP\vert}$ such that the following diagram commutes:
\begin{equation}\label{eq:bundles2}
\xymatrix{
T_\lambda\ar[rr]^L\ar[rd]_\pi & & \reals^{\vert\cP\vert}\times(\reals^+)^{\vert\cP\vert}\ar[ld]\\
 & (\reals^+)^{\vert\cP\vert} & }
 \end{equation}
Here the unlabeled arrow is the projection on the first factor.

Combining \eqref{eq:bundles} and \eqref{eq:bundles2} we obtain a trivialization of the principal bundle $\ell_{\cP}:\cT(X)\to(\reals^+)^{\vert\cP\vert}$, i.e. Fenchel-Nielsen coordinates, which differ from the shearing coordinates $s_\lambda$ by a linear map. Since by Theorem \ref{thm:shear-convex} the length functions are convex with respect to the shearing coordinates and convexity is preserved by linear maps, the same result holds for this choice of Fenchel-Nielsen coordinates. This concludes the proof of Theorem \ref{thm:FN-convex}.\qed{thm:FN-convex}

%
%
%
%

\section{Length functions on $\cT(A)$}\label{sec:annulus}
Let $A$ be a complete hyperbolic annulus and assume that the core curve $\alpha$ is isotopic to a geodesic and let $\cT(A)$ be the Teichm\"uller space of $A$. Observe that according to the definition above, $\cT(A)$ is infinite dimensional. Fix also an ideal triangulation $\tilde\lambda$ of $A$ and let $\lambda$ be the sublamination of $\tilde\lambda$ obtained by deleting all leaves of $\lambda$ which are disjoint from the core geodesic of $A$. Assume that $\lambda$ contains only finitely many non-isolated leaves. In this section we construct a map 
$$s_\lambda:\cT(A)\to\reals^{\vert\lambda\vert}$$
and a convex function
$$\bar{\ell}: \reals^{|\lambda|} \longrightarrow (0, \infty]$$
such that $\ell = \bar{\ell} \circ s_\lambda$. Here 
$$\ell : \cT(A) \longrightarrow (0,\infty)$$
is the function which assigns to each hyperbolic annulus the length of its core curve.

\subsection{A treatise on wedges}
An {\em ideal wedge} is the region bounded by two asymptotic geodesics in $\htwo$. If $W$ is an ideal wedge in $\htwo$ or $A$ then there is a unique ideal triangle $\Delta(W)$ two of whose boundary components are the boundary components of $W$. In particular, the two boundary edges $\gamma_1$ and $\gamma_2$ of the ideal wedge $W$ have well-determined {\em midpoints}. Orienting both boundary edges we parametrize $\gamma_1$ and $\gamma_2$ by $\reals$ via the signed distance from their. To $(x,y)\in\reals^2$ we associate the pair of points in $\bdry W$ corresponding to $x$ in $\gamma_1$ and $y$ in $\gamma_2$ respectively; let $d(x,y)$ be he distance in $\htwo$ of these two points in $\bdry W$.


\begin{figure}[h]
\leavevmode \SetLabels
\L (.26*.68) $(x,y)$\\
\L (.67*.62) $y$\\
\L (.58*.40) $x$\\
\endSetLabels
\par
\begin{center}
\AffixLabels{\centerline{\epsfig{file =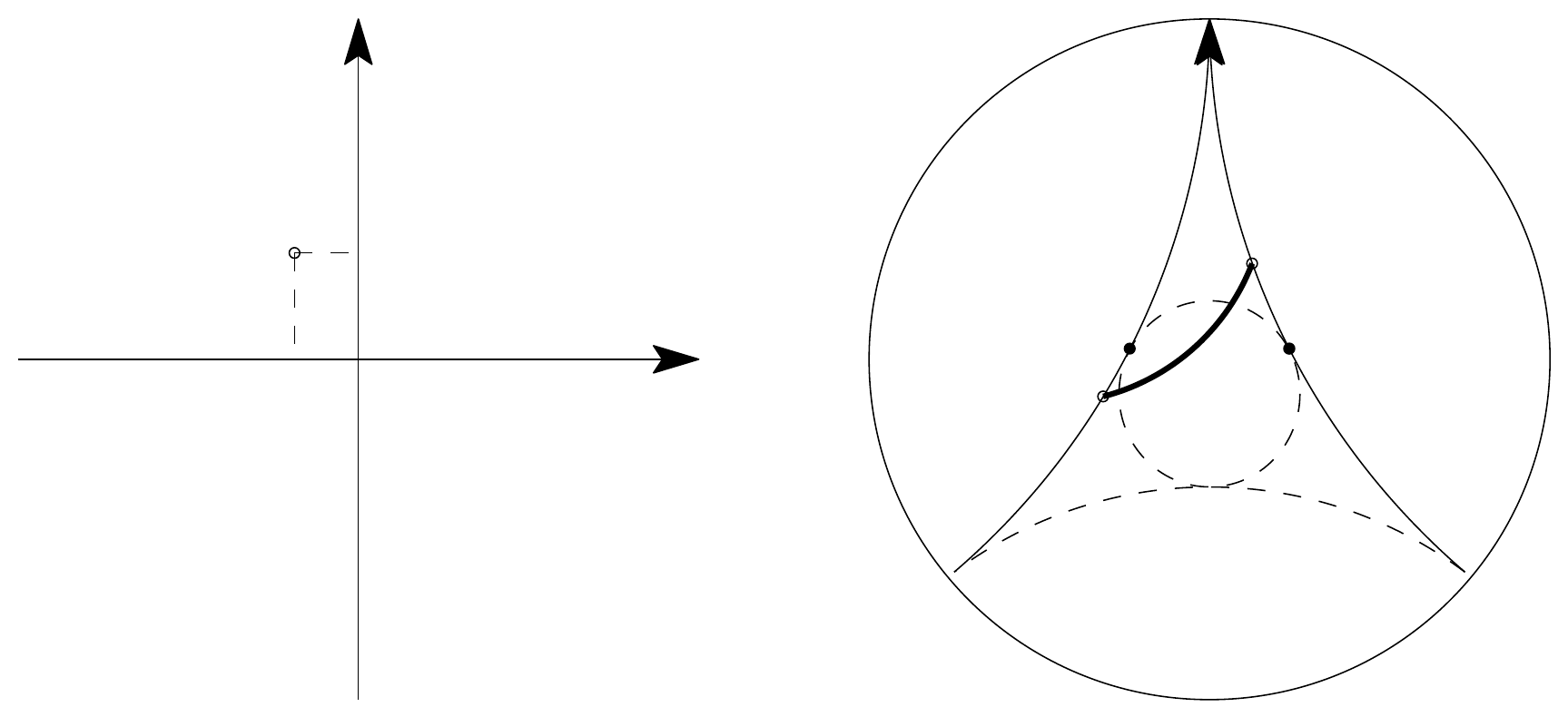, width=8cm, angle= 0}}}
\end{center}
\caption{The length of the bold segment is by definition $d(x,y)$.}\label{fig4}
\end{figure}

We can describe the function
$$d: \reals^2 \longrightarrow (0,\infty)$$
more explicitly using the upper half space model of $\htwo$. Namely, define $d(x,y)$ to be the distance in $\htwo$ between the points $\imath e^x$ and $1 + \imath e^y$.

\begin{lemma}\label{lem:basic-convexity}
The function $d$ is strictly convex and
$$|x-y| \leq d(x,y) \leq |x-y| + \frac{1}{\max\{e^x, e^y\}}.$$
\end{lemma}

Lemma \ref{lem:basic-convexity} can be proved via a simple computation but it should be remarked that in more intrinsic terms it just follows from the fact that $\htwo$ is negatively curved and hence that the distance function is convex.\qed{lem:basic-convexity}

An  injective isometric immersion of a wedge into a hyperbolic annulus $A$ will also be called an ideal wedge; all the definitions above carry over without difficulties.
\medskip

An  {\em ideal {\wln}} is a geodesic lamination $\lambda$ on $A$ whose complement is a disjoint union of ideal wedges, and such that $\lambda$ contains only finitely many non-isolated leaves. For the sake of concreteness we will also assume that every leaf of $\lambda$ which is isolated to one side is actually isolated. 

\begin{lemma}\label{lem:lift-wln}
Let $S$ be a hyperbolic surface and $\lambda$ an ideal triangulation of $S$. Let also $A$ be an annulus and $\pi:A\to S$ a covering. The set of leaves of $\pi^{-1}(\lambda)$ which intersect every curve in $A$ homotopic to the core curve is an ideal {\wln}.\qed{lem:lift-wln}
\end{lemma}

On a finite area surface it is not possible to consistently orient an ideal triangulation. On an annulus, an ideal {\wln} can be consistently oriented and this will will be important in the work below. To do so we fix an orientation on $A$ and of its core geodesic $\alpha$. Observe that all the leaves of an ideal {\wln} $\lambda$ of $A$ intersect $\alpha$ exactly once. We then orient every geodesic in $\lambda$ so that if $\gamma$ is a geodesic in $\lambda$ the orientation of $A$ at $\alpha \cap \gamma$ is given by the ordered pair of the oriented tangent vectors to $\alpha$ and $\gamma$. 

\begin{figure}[h]
         \centering
         \includegraphics[width=6cm]{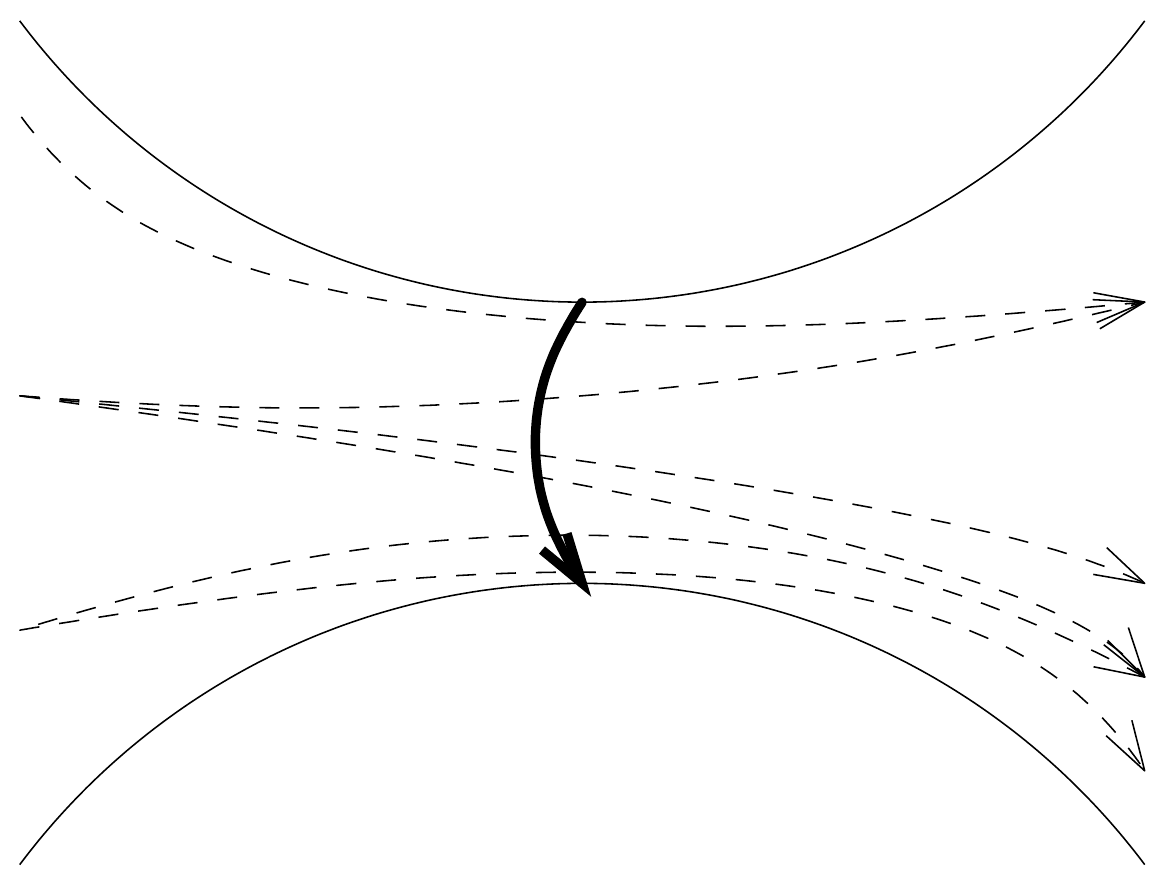}
         \caption{An ideal {{\wln}} of the annulus $A$ with the corresponding orientations induced by the bold printed core geodesic.}\label{fig3}
\end{figure}

Observe that if $W$ is an ideal wedge of in $A\setminus\lambda$, then the orientations of the boundary edges of $W$ determined by $W$ and the orientations determined by $\lambda$ agree if and only if the vertex of $W$ is to the left of the core curve $\alpha$.

\subsection{The shearing map}
Recall that for every ideal wedge $W$ in $\htwo$ or $A$ we have associated a unique ideal triangle $\Delta(W)$. If $W^a$ and $W^b$ are disjoint wedges in $\htwo$ and $\gamma$ is a geodesic in $\htwo$ separating $W^a$ from $W^b$ and asymptotic to some boundary edge of $W^a$ and some boundary edge $W^b$ we define
$$s(W^a, W^b, \gamma) = s(\Delta(W^a), \Delta(W^b), \gamma).$$
Clearly, all claims of Lemma \ref{lem:facts} hold in this setting as well.

An ideal {\wln} of $A$ yields now a map
$$s_\lambda : \cT(A) \longrightarrow \reals^{|\lambda|}$$
in the same way as it did for an ideal triangulation on a surface except we replace ideal triangles with ideal wedges $W^-_\gamma$ and $W^+_\gamma$. We remark that $s_\lambda$ does not yield coordinates of $\cT(A)$. We refer to the map $s_\lambda$ as the {\em shearing map}.

\subsection{The function $L$}
Let $\lambda$ be an ideal {\wln} of $A$, $\lambda_0\subset\lambda$ the collection of isolated leaves of $\lambda$ and $\cW$ the set of ideal wedges in $A \setminus \lambda$. Recall that $\lambda\setminus\lambda_0$ consists by definition of only finitely many leaves and that we are also working under the additional assumption that the leaves in $\lambda\setminus\lambda_0$ are non-isolated on both sides.

For each $W \in \cW$ let $\gamma^-_W$ be the left boundary edge and $\gamma^+_W$ be the right boundary edge where left and right is defined with respect to the orientation of the core curve $\alpha$. Observe that both $\gamma^+_W$ and $\gamma^-_W$ belong to $\lambda$. If the two boundary components $\gamma^+_W$ and $\gamma^-_W$ of $W$ are asymptotic to the left $\alpha$ then define $d_W(x,y) = d(x,y)$ for $(x,y)\in\reals^2$. If they are asymptotic to the right of $\alpha$ then define $d_W(x,y) = d(-x,-y)$. Here left and right are defined with respect to the orientation of $\alpha$.

Define a function $L: \reals^{|\lambda|} \times \reals^{|\lambda_0|} \longrightarrow (0, \infty]$ by
$$L(\x, \y) = \sum_{W \in \cW} d_W(\gamma^-_W(\y), \gamma^+_W(\y) + \gamma^+_W(\x)).$$
Here $\gamma(\x)$ is the $\gamma$-coordinate of $\x\in\reals^\lambda$. Note that $L$ does not depend on the coordinates of the non-isolated leaves in $\lambda$.
\medskip

\noindent {\bf Meaning of $L$:}
The definition of $L$ is a bit obscure. Before moving on we explain its meaning in the particular case that $\lambda$ is a finite {\wln}. In particular $\lambda_0=\lambda$ and the sum in the definition of $L$ is finite. Given a point $(B,f)\in\cT(A)$, its image $\x=s_\lambda(B,f)$ under the shearing map encodes how the different wedges are glued in $X$. The element $\y\in\reals^{\lambda_0}$ picks a point in each leaf of the (finite) {\wln} $\lambda$. In particular, $\y$ picks in every wedge $W$ two different points in the boundary. Let $\alpha_W$ be the geodesic segment in $W$ joining these two points. The juxtaposition of all the segments $\alpha_W$ is a closed loop in $(B,f)$ homotopic to the core curve; $L(\x,\y)$ is the length of this loop.

\subsection{Piecewise geodesic segments}
Our next aim is to generalize the preceding discussion on the meaning of the function $L$. Assume now that $\lambda$ is a general {\wln} and recall that $\lambda_0$ is the collection of isolated leaves of $\lambda$. 

For a marked hyperbolic annulus  $(B,f) \in \cT(A)$ let $\cP(B,f)$ be the set of closed continuous loops $\beta$ on $B$ homotopic to $f(\alpha)$ and such that the intersection of $\beta$ with each wedge in $B \backslash \bar{f}(\lambda)$ is a geodesic segment. In particular, any $\beta\in\cP(B,f)$ intersects each geodesic in $\bar{f}(\lambda)$ exactly once. We let 
$$p_\lambda: \cP(B,f) \longrightarrow \reals^{|\lambda_0|}$$ 
be the map defined by the property that $\gamma(p_\lambda(\beta))$ is the signed distance between $m_{\bar{f}(\gamma)}$ and the intersection of $\beta$ with $\bar{f}(\gamma)$ for all $\gamma \in \lambda_0$. As always, the sign is determined by the orientation induced by the core geodesic $\alpha$.

The map $p_\lambda$ is clearly injective. However, if $\lambda_0\neq\lambda$ it is easy to see that it is not surjective. Our next goal is to determine the image of $p_\lambda$. In order to do so we define two quantities $a^-_\gamma(\x, \y)$ and $a^+_\gamma(\x,\y)$ for every non-isolated leaf $\gamma$; recall that we are assuming that leaves which are isolated to one side are actually isolated. We first define $a^-_\gamma$. Label the geodesics between $\gamma$ and $W^-_\gamma$, the (arbitrarily) chosen wedge to the left of $\gamma$, by $\gamma^-_0=\gamma^+_{W^-_\gamma}, \gamma^-_1, \dots$ so that they are positively ordered with respect to $\alpha$. Define 
\begin{equation}\label{left-limit}
a^-_\gamma(\x,\y) = \underset{n \to \infty}{\lim}\left(\gamma^-_n(\y) + \sum_{i=0}^n \gamma^-_i(\x)\right)
\end{equation}
if the limit exists. The definition of $a^+_\gamma(\x,\y)$ is very similar with small changes. Label the geodesics between $W^+_\gamma$ and $\gamma$ by $\gamma^+_0=\gamma^-_{W^+_\gamma}, \gamma^+_1, \dots$ so that they are negatively ordered with respect to $\gamma$ and define
\begin{equation}\label{right-limit}
a^+_\gamma(\x,\y) = \underset{n \to \infty}{\lim}\left(\gamma^+_n(\y) - \sum_{i=0}^{n-1} \gamma^+_i(\x)\right)
\end{equation}
if the limit exist and $a^+_\gamma(x,y)=\infty$ otherwise. To understand the meaning of the sums in the definition compare with the last claim in Lemma \ref{lem:facts}.

Given $\x\in\reals^{\vert\lambda\vert}$ define $P_{\x} \subset \reals^{|\lambda_0|}$ to be the set
$$P_{\x} = \left\{\y \in \reals^{|\lambda_0|} {\Big|} 
\begin{array}{c}
\hbox{the limits \eqref{left-limit} and \eqref{right-limit} exist, and} \\
 a^-_\gamma(\x, \y) = a^+_\gamma(\x, \y) + \gamma(\x)< \infty \mbox{ for all } \gamma \in \lambda \backslash \lambda_0
\end{array}
\right\}.$$
We prove:

\begin{lemma}\label{curvefamily}
Let $(B,f)$ be a marked hyperbolic annulus in $\cT(A)$. If $\x = s_\lambda(B,f)$ then the $p_\lambda$-image of $\cP(B,f)$ is $P_\x$ and
$$\length(\beta)  = L(\x, p_\lambda(\beta)).$$
\end{lemma}

\noindent {\bf Proof.} Let $\gamma$ be a geodesic in $\lambda \backslash \lambda_0$. As above we label the geodesics on the left of $\gamma$ by $\gamma^-_0, \gamma^-_1, \dots$ and on the right by $\gamma^+_0, \gamma^+_1, \dots$. If $\beta$ is a curve in $\cP(B,f)$ let $p^\pm_i$ be the point of intersection of $\beta$ with $\bar{f}(\gamma^\pm_i)$. Note that the sequence $p^-_i$ and $p^+_i$ will limit to the same point $p$ on $\gamma$ since $\beta$ is a continuous path. Therefore both $h_{\bar{f}(\gamma^-_i), \bar{f}(\gamma)}(p^-_i)$ and $h_{\bar{f}(\gamma^+_i), \bar{f}(\gamma)}(p^+_i)$ will limit to $p$. Let $d^-_i$ be the signed distance between $h_{\bar{f}(\gamma^-_0), \bar{f}(\gamma)}(m_{\bar{f}(\gamma^-_0)})$ and $h_{\bar{f}(\gamma^-_i), \bar{f}(\gamma)}(p^-_i)$ and similarly define $d^+_i$. Then $d^-_i$ will limit to the signed distance between $h_{\bar{f}(\gamma^-_0), \bar{f}(\gamma)}(m_{\bar{f}(\gamma^-_0)})$ and $p$ which we label $d^-$. Similarly $d^+$, the limit of $d^+_i$, will be the signed distance between $h_{\bar{f}(\gamma^+_0), \bar{f}(\gamma)}(m_{\bar{f}(\gamma^+_0)})$ and $p$. Since the signed distance between $h_{\bar{f}(\gamma^-_0), \bar{f}(\gamma)}(m_{\bar{f}(\gamma^-_0)})$ and $h_{\bar{f}(\gamma^+_0), \bar{f}(\gamma)}(m_{\bar{f}(\gamma^+_0)})$ is $\gamma(\x)$ we have
\begin{equation}\label{limitsum}
d^- - d^+ = \gamma(\x).
\end{equation}
Since the $h$-maps are isometries and $h_{\bar{f}(\gamma^-_i), \bar{f}(\gamma)} \circ h_{\bar{f}(\gamma^-_0), \bar{f}(\gamma^-_i)} = h_{\bar{f}(\gamma^-_0), \bar{f}(\gamma)}$ we have that $d^-_i$ is also equal to the signed distance between $h_{\bar{f}(\gamma^-_0), \bar{f}(\gamma^-_i)}(m_{\bar{f}(\gamma^-_0)})$ and $p^-_i$. By definition $\gamma^-_n(p_\lambda(\beta))$ is the signed distance between $m_{\bar{f}(\gamma^-_n)}$ and $p^-_n$ so we have
$$d^-_n =  \gamma^-_n(p_\lambda(\beta)) + \sum_{i=0}^n \gamma^-_i(\x)$$
by the last claim of Lemma \ref{lem:facts}. Therefore $d^- = a^-_\gamma(\x, p_\lambda(\beta))$ and similarly $d^+ = a^+_\gamma(\x, p_\lambda(\beta))$. Rearranging \eqref{limitsum} we have
$$a^-_\gamma(\x, p_\lambda(\beta)) = a^+_\gamma(\x, p_\lambda(\beta)) + \gamma(\x)$$
and therefore $p_\lambda(\beta)$ is in $P_\x$.

For every $\y \in P_\x$ we need to build a curve $\beta$ on $B$ such that $p_\lambda(\beta) = \y$. Given a geodesic $\alpha$ in $\lambda_0$ let $p_\alpha$ be the point on $\bar{f}(\alpha)$ whose signed distance from $m_{\bar{f}(\alpha)}$ is $\alpha(\y)$. Then on each of the finitely many components of $B \backslash \bar{f}(\lambda \backslash \lambda_0)$ there is an arc that intersects each $\alpha$ in $\lambda_0$ at $p_\alpha$ and is a geodesic in each wedge. For these arcs to complete to a simple closed curve in $\cP(B,f)$ we need that for each geodesic $\gamma$ in $\lambda \backslash \lambda_0$ the limit of the points $p_{\gamma^-_i}$ is equal to the limit of $p_{\gamma^+_i}$. From the first paragraph of this proof we see that this holds when $\y$ is in $P_\x$.

Finally we note that the length of $\beta$ is the sum of the lengths of the restriction of $\beta$ to each ideal wedge. This is exactly the sum $L(\x, p_\lambda(\beta))$. \qed{curvefamily}

Before moving on to more interesting topics, we observe:

\begin{lemma}\label{lem:subspace}
The set $\bb P=\{(\x,\y)\in\reals^{|\lambda|}\times\reals^{|\lambda_0|}\vert \y\in P_{\x}\}$ is a linear subspace and the projection $\bb P\to\reals^{\vert\lambda\vert}$ is surjective.\qed{lem:subspace}
\end{lemma}

\subsection{Convexity}
We restrict from now on the function $L$ to the linear subspace $\bb P$. The function $L$ is an infinite sum of convex functions. This sum will not be finite everywhere so we need an extended notion of convexity. If $V$ is vector space and $f: V \longrightarrow (-\infty, \infty]$ is a function then $f$ is {\em convex} if
$$f\left(\frac{\x_0 + \x_1}{2} \right) \leq \frac{f(\x_0) + f(\x_1)}{2}$$
whenever $f(\x_0)$ and $f(\x_1)$ are finite. The function is {\em strictly convex} if the inequality is always strict. We need the following general lemma.
\begin{lemma}\label{infconvexity}
Let $V_0$ and $V_1$ be vector spaces and $P$ a subspace of $V_0 \times V_1$ whose projection onto $V_0$ is onto. Let $P_\x = \{ \y \in V_1 | (\x,\y) \in P\}$.
If $F: V_0 \times V_1 \longrightarrow (0, \infty]$ is convex then
$$f(\x) = \underset{\y \in P_\x}{\inf} F(\x,\y)$$
is convex. If
$$f(\x) = \underset{\y \in P_\x}{\min} F(\x,\y)$$
and $F$ is strictly convex on $P$ then $f$ is strictly convex.
\end{lemma}

\noindent {\bf Proof.} Let $\x_0$ and $\x_1$ be points in $V_0$. We can assume that both $f(x_0)$ and $f(x_1)$ are finite for otherwise the lemma is trivial. For any $\epsilon>0$ we can choose a $\y_0 \in P_{\x_0}$ and $\y_1 \in P_{\x_1}$ such that $f(\x_i) > F(\x_i, \y_i) - \epsilon$. Since $P$ is a subspace we have that $(\y_0 + \y_1)/2$ is in $P_{(\x_0+ \x_1)/2}$ and therefore
\begin{eqnarray*}
f\left(\frac{\x_0 + \x_1}{2}\right) &\leq & F\left(\frac{\x_0+\x_1}{2},  \frac{\y_0 + \y_1}{2}\right) \\
& \leq & \frac{F(\x_0,\y_0) + F(\x_1, \y_1)}{2} \\
& \leq & \frac{f(\x_0) + f(\x_1)}{2} + \epsilon.
\end{eqnarray*}
Since the $\epsilon$ is arbitrary we have that
$$f\left(\frac{\x_0 + \x_1}{2}\right) \leq  \frac{f(\x_0) + f(\x_1)}{2}$$
and $f$ is convex.

If $f$ is defined by a minimum then there are $\y_0$ and $\y_1$ such that $f(\x_i) = F(\x_i, \y_i)$. If $F$ is also strictly convex we have
\begin{eqnarray*}
f\left(\frac{\x_0 + \x_1}{2}\right) &\leq & F\left(\frac{\x_0+\x_1}{2},  \frac{\y_0 + \y_1}{2}\right) \\
& < & \frac{F(\x_0,\y_0) + F(\x_1, \y_1)}{2} \\
& = & \frac{f(\x_0) + f(\x_1)}{2}
\end{eqnarray*}
which implies that $f$ is strictly convex. \qed{infconvexity}

\begin{theorem}\label{annulusconvex}
Let $\bar{\ell}: \reals^{|\lambda|} \longrightarrow (0, \infty]$ be defined by
$$\bar{\ell}(\x) = \underset{\y \in P_\x}{\min} L(\x,\y).$$
Then $\bar{\ell}$ is strictly convex and
$$\ell = \bar\ell \circ s_\lambda.$$
\end{theorem}

\noindent {\bf Proof.} 
To see that the minimum defining $\bar{\ell}$ exists we use Lemma \ref{curvefamily}. At the same time we will see that $\ell$ is the composition of $s_\lambda$ and $\bar{\ell}$. The set $\cP(B,f)$ contains the geodesic representative in the homotopy class of $f(\alpha)$ so
\begin{eqnarray*}
\ell(B,f) & = &\underset{\beta \in \cP(B,f)}{\min} \length(\beta) \\
& = & \underset{\beta \in \cP(B,f)}{\min} L(s_\lambda(B,f), p_\lambda(\beta)) \\
& = & \underset{\y \in P_\x}{\min} L(\x, \y) \\
& = & \bar{\ell} \circ s_\lambda(B,f)
\end{eqnarray*}
where second two inequalities follow from Lemma \ref{curvefamily} and $\x = s_\lambda(B,f)$.

Let $W$ be a wedge in $\cW$.
Since the function $d$ is strictly convex the function $(\x,\y) \mapsto d_W(\gamma^-_W(\y), \gamma^+_W(\y) + \gamma^+_W(\x))$ is convex on $\reals^{|\lambda|} \times \reals^{|\lambda_0|}$ and is strictly convex along any line through distinct points $(\x_0, \y_0)$ and $(\x_1, \y_1)$ were either $\gamma^+_W(\x_0) \neq  \gamma^+_W(\x_1)$, $\gamma^+_W(\y_0) \neq \gamma^+_W(\y_0)$ or $\gamma^-_W(\y_0)\neq \gamma^-_W(\y_1)$. Therefore if $L$ is not strictly convex along this line we must have $\beta(\x_0) = \beta(\x_1)$ and $\beta(\y_0) = \beta(\y_1)$ for all $\beta \in \lambda_0$. For all closed curves $\gamma$ in $\lambda$ we then have $a^\pm_\gamma(\x_0, \y_0)  = a^\pm_\gamma(\x_1, \y_1)$ and therefore if $\y_0 \in P_{\x_0}$ and $\y_1 \in P_{\x_1}$ we also have $\gamma(\x_0) = \gamma(\x_1)$. Therefore if $L$ is not strictly convex along the line it can't lie in $\bb P$; equivalently $L$ is strictly convex on $\bb P$. Lemma \ref{infconvexity} then implies that $\bar{\ell}$ is strictly convex.
\qed{annulusconvex}

%
%
%
%

\section{Proof of Proposition \ref{lengthfactors}}\label{sec:final}

We now return to a finite area hyperbolic surface $X$ with an ideal triangulation $\lambda$ and recall the proposition.
\medskip

\noindent {\bf Proposition \ref{lengthfactors}} 
{\em Let $\alpha$ be a closed curve on $X$ and $\ell_\alpha : \cT(X) \longrightarrow \reals^+$ its length function. Then there is a convex function $\bar{\ell}_\alpha: \reals^{|\lambda|} \longrightarrow \bar{\reals}$ such that $\ell_\alpha = \bar{\ell}_{\alpha} \circ s_\lambda$. Furthermore the $\bar{\ell}_{\alpha}$-image of $T_\lambda$ is in $\reals^+$ and if $\alpha$ intersects every leaf of $\lambda$ then $\bar{\ell}_\alpha$ is strictly convex.}
\medskip

Let $\pi_\alpha: A_\alpha \longrightarrow X$ be the annular covering such that $\alpha$ lifts to $A_\alpha$. If $(Y,f)$ is a pair in $\cT(X)$ then $Y_{f(\alpha)}$ is the annulus cover of $Y$ such that $f(\alpha)$ lifts to $Y_{f(\alpha)}$. The map $f:X\to Y$ lifts to a quasi-conformal homeomorphism $f_\alpha: X_\alpha \longrightarrow Y_{f(\alpha)}$. This induces a map $\Pi_\alpha : \cT(X) \longrightarrow \cT(X_\alpha)$. Let $\ell_\alpha: \cT(X) \longrightarrow (0,\infty)$ be the geodesic length function for $\alpha$ and let $\ell: \cT(X_\alpha) \longrightarrow (0,\infty)$ be the geodesic length function for the core curve of the annulus. Then $\ell_\alpha = \ell \circ \Pi_\alpha$.

Let $\lambda_\alpha$ be the set of geodesics in $\pi^{-1}_\alpha(\lambda)$ that intersect the core curve $\alpha$. Then $\lambda_\alpha$ is an ideal {\wln} of $X_\alpha$ by Lemma \ref{lem:lift-wln}. The covering map $\pi_\alpha$ defines a map $\bar{\pi}_\alpha$ from the leaves of $\lambda_\alpha$ to the leaves of $\lambda$ which induces a map $\bar{\Pi}_\alpha : \reals^{|\lambda|} \longrightarrow \reals^{|\lambda_\alpha|}$.

\begin{prop}\label{coordinatescommute}
We can choose the coordinates for $\cT(A)$ and $\cT(S)$ such that $s_{\lambda_\alpha} \circ \Pi_\alpha = \bar{\Pi}_\alpha \circ s_\lambda$.
\end{prop}

\noindent{\bf Proof.} The covering map $\pi:\htwo \longrightarrow X$ factors through $\pi_\alpha$. In particular $X$ and $A_\alpha$ have the same universal cover and the pre-image $\tilde{\lambda}_\alpha$ of $\lambda_\alpha$ is a subset of the pre-image $\tilde{\lambda}$ of $\lambda$. Assume we have chosen the coordinate map $s_\lambda$.
Let $\gamma'$ be a geodesic in $\lambda_\alpha$ and let $\gamma = \pi_\alpha(\gamma')$ be its image in $\lambda$. Let $\tilde{\gamma}'$ be a component of the pre-image of $\gamma'$ in $\htwo$. In our choice of coordinate map $s_\lambda$ we have chosen a pre-image $\tilde\gamma$ in $\htwo$ of $\gamma$ and ideal triangles $\Delta^a_\gamma$ and $\Delta^b_\gamma$.
 Then there is a deck transformation $g \in \pi_1(X)$ such that $g(\tilde{\gamma}) = \tilde{\gamma}'$. Since $\tilde{\gamma}$ separates $\Delta^a_\gamma$ from $\Delta^b_\gamma$ we have that $\tilde{\gamma}'$ separates $g(\Delta^a_\gamma)$ from $g(\Delta^b_\gamma)$. These ideal triangles will also be asymptotic to $\tilde{\gamma}'$ so they must have exactly two edges that intersect $\tilde{\alpha}$, the pre-image of $\alpha$ in $\htwo$ and these edges will be in $\tilde{\lambda}_\alpha$. Exactly one of these ideal triangles, say $g(\Delta^a_\gamma)$, will be on the left $\tilde{\gamma}'$. Since it has two edges in $\tilde{\lambda}_\alpha$ there will be a wedge $W$ in $\htwo \backslash \tilde{\lambda}_\alpha$ with $\Delta(W) = g(\Delta^a_\gamma)$. We let $W^-_\gamma= W$. Similarly the ideal triangle $g(\Delta^b_\gamma)$ determines a wedge in $\htwo\backslash \tilde{\lambda}_\alpha$ which we define to be $W^+_\gamma$. If $g(\Delta^b_\gamma)$ is on the left we reverse the labels. It is then straightforward to check that the lemma holds. \qed{coordinatescommute}

Let $\bar{\ell} : \reals^{|\lambda_\alpha|} \longrightarrow [0,\infty]$ be the factorization of $\ell$ given in the previous subsection so that $\ell = \bar{\ell} \circ s_{\lambda_\alpha}$. Define $\bar{\ell}_\alpha : \reals^{|\lambda|} \longrightarrow [0, \infty]$ by $\bar{\ell}_\alpha = \bar{\ell} \circ \bar{\Pi}_\alpha$. We then have the following commutative diagram:

$$\xymatrix{
\cT(X)\ar[rr]^{\sigma_\lambda}\ar[dd]^{\Pi_\alpha}\ar[rd]^{\ell_\alpha} & & \reals^{|\lambda|}\ar[dd]^{\bar\Pi_\alpha}\ar[ld]^{\bar\ell_\alpha} \\
& [0,\infty] & \\
\cT(X_\alpha)\ar[rr]_{\sigma_{\lambda_\alpha}}\ar[ru]_\ell & &\reals^{|\lambda_\alpha|}\ar[lu]^{\bar\ell}
 }
 $$
The first claim of Proposition \ref{lengthfactors} follows directly from the commutativity of the triangle on the top. The convexity of $\bar\ell_\alpha$ follows from the strict convexity of $\bar\ell$ and the fact that $\bar\Pi_\alpha$ is linear. In the case that $\alpha$ intersects every leaf of $\lambda$, the map $\bar\Pi_\alpha$ is also injective and hence $\bar\ell_\alpha$ is convex as claimed. The following lemma now concludes the proof of Proposition \ref{lengthfactors}:

\begin{lemma}\label{finite}
The function $\bar{\ell}_\alpha$ is finite on $T_\lambda$.
\end{lemma}

If we knew at this stage that $s_\lambda$ is a homeomorphism onto its image $T_\lambda$, the claim of Lemma \ref{finite} would follow also from the commutativity of the diagram above and the fact that $\ell_\alpha$ is obviously finite. Unfortunately, we derived this fact from Proposition \ref{lengthfactors}.
\medskip

\noindent {\bf Proof.} Let $\x$ be in $T_\lambda$ and let $\x_\alpha = \bar{\Pi}_\alpha(\x)$. We will choose a $\y \in \reals^{|(\lambda_\alpha)_0|}$ such that $\y \in \cP_{\x_\alpha}$ and $L(\x_\alpha, \y) < \infty$. Let $\gamma_\alpha$ be a non-isolated leaf of $\lambda_\alpha$. Then $\bar{\pi}_\alpha(\gamma_\alpha)$ is a closed leaf $\gamma$ of $\lambda$. Let $\gamma^+_i$ be the leaves on the right of $\gamma_\alpha$ labeled as the were in the previous section. Assume that the right side of $\gamma_\alpha$ maps to the $a$-side of $\gamma$. Then $\bar{\pi}_\alpha(\gamma^+_i) = \bar{\pi}_\alpha(\gamma^+_{i+k})$
where
$$k = \sum_{\beta \in \lambda} n^a_\gamma(\beta).$$
Therefore $\gamma^+_i(\x_\alpha) = \gamma^+_{i+k}(\x_\alpha)$ and there is a constant $c$ such that
$$c =  \sum_{i=j}^{j+k-1} \gamma^+_i(\x_\alpha)$$
for all $j\geq 0$. It follows that
$$\sum_{i = 0}^{nk +j} \gamma^+_i(\x_\alpha) = nc + \sum_{i=0}^{j}(\x_\alpha).$$

Choose the $\gamma^+_i$-coordinates of $\y$ such that
$$\gamma^+_n(\y) =  \sum_{i=0}^{n-1} \gamma^+_i(\x_\alpha) = \gamma^+_{n-1}(\y) + \gamma^+_{n-1}(\x_\alpha).$$
Let $\cW^+_{\gamma_\alpha}$ be the set of wedges in $\cW$ whose boundary edges are in $\gamma^+_i$ and $\gamma_i^+$ for some $i\ge 0$. Then
\begin{eqnarray*}
\sum_{W \in \cW^+_{\gamma_\alpha}} d_W(\gamma^-_W(\y), \gamma^+_W(\y) + \gamma^+_W(\x_\alpha)) &=& \sum_{j=0}^{k-1} \sum_{i=0}^\infty d\left(\gamma^+_{j +ki+1}(\y), \gamma^+_{{j+ki}}(\y) + \gamma^+_{{j+ki}}(\x_\alpha)\right) \\
& = & \sum_{j=0}^{k-1} \sum_{i=0}^\infty d\left(ic + \gamma_{j}^+(\x_\alpha), ic + \gamma^+_{j}(\x_\alpha)\right) \\
& \leq & \sum_{j=0}^{k-1} \sum^\infty_{i=0} e^{-ic - \gamma^+_{j}(\x_\alpha)}\\
&<& \infty
\end{eqnarray*}
where the inequality on the third line follows from Lemma \ref{lem:basic-convexity}.

We now take the leaves $\gamma^-_i$ on the left of $\gamma_\alpha$. We similarly choose the $\gamma^-_i$-coordinates but so that we have $a^+_{\gamma_\alpha}(\x_\alpha, \y) + \gamma_\alpha(\x_\alpha) = a^-_{\gamma_\alpha}(\x_\alpha, \y)$ we define
$$\gamma^-_n(\y)  = \gamma_\alpha(\x_\alpha)-\sum_{i=0}^n \gamma^-_i(\x_\alpha).$$
Similar to above we define $\cW^-_{\gamma_\alpha}$ to be the set of wedges in $\cW$ whose boundary edges are $\gamma^-_i$ and $\gamma^-_{i+1}$ for some $i \ge 0$ we similarly find that
$$\sum_{W \in \cW^-_{\gamma_\alpha}} d_W(\gamma^-_W(\y), \gamma^+_W(\y) + \gamma^+_W(\x_\alpha)) < \infty.$$

We repeat the above construction for each of the finitely many non-isolated leaves of $\lambda_\alpha$.
Note that each wedge $W \in \cW$ is asymptotic to at most one of the non-isolated leaves $\gamma_\alpha$ and therefore is in at most one of the sets $\cW^+_{\gamma_\alpha}$ or $\cW^-_{\gamma_\alpha}$ that contains $W$ and all but finitely many of the $W$ are in one of these sets.  This defines $\y$ for all but finitely many coordinates. The remaining we choose freely.  We can then split the infinite sum $L(\x_\alpha, \y)$ into finitely many sums each of which is finite. Therefore $\bar{\ell}_\alpha(\x) \leq L(\x_\alpha, \y) < \infty$. \qed{finite}

This concludes the proof of Proposition \ref{lengthfactors} and with it the whole paper.

\bigskip

\noindent K. Bromberg\newline
\noindent Department of Mathematics, University of Utah
\newline \noindent
\texttt{bromberg@math.utah.edu}
\bigskip

\noindent M. Bestvina
\newline \noindent Department of Mathematics, University of Utah
\newline \noindent
\texttt{bestvina@math.utah.edu}
\bigskip

\noindent K. Fujiwara
\newline\noindent Graduate School of Information Sciences, Tohoku University 
\newline \noindent
\texttt{fujiwara@math.is.tohoku.ac.jp}
\bigskip

\noindent J. Souto
\newline \noindent Department of Mathematics, University of Michigan
\newline \noindent
\texttt{jsouto@umich.edu}

 \end{document}

%% file: macros.tex
\newcommand{\bb}{\mathbb}

\newcommand{\reals}{{\bb R}}

\newcommand{\htwo}{{{\bb H}^2}}

\newcommand{\widemargins}{
\setlength{\textwidth}{5.8in}
\setlength{\oddsidemargin}{0.25in}
\setlength{\evensidemargin}{0.25in}
}

\newcommand{\qed}[1]{\nopagebreak[4]{\tiny \hfill
\fbox{\ref{#1}} \linebreak }\pagebreak[2]}
\newcommand{\bdry}{\partial}

\newcommand{\del}{\partial}

\newcommand{\area}{\operatorname{area}}

\newcommand{\length}{\operatorname{length}}

\newtheorem{theorem}{Theorem}[section]
\newtheorem{prop}[theorem]{Proposition}
\newtheorem{lemma}[theorem]{Lemma}
\newtheorem{cor}[theorem]{Corollary}

\newcommand{\cP}{{\cal P}}

\newcommand{\cT}{{\cal T}}

\newcommand{\cW}{{\cal W}}
